%% file: paper24.tex
\documentclass[12pt]{amsart}
\usepackage{amscd,amssymb,a4wide,color,graphics,pstricks}

\usepackage{mathrsfs}
\input xy
\xyoption{all}

\newtheorem{theorem}{Theorem}[section]
\newtheorem{lemma}[theorem]{Lemma}
\newtheorem{proposition}[theorem]{Proposition}
\newtheorem{corollary}[theorem]{Corollary}
\newtheorem{claim}[theorem]{Claim}

\theoremstyle{definition}
\newtheorem{definition}[theorem]{Definition}

\newtheorem{observation}[theorem]{Observation}

\theoremstyle{remark}
\newtheorem{remark}[theorem]{Remark}

\newcommand{\NN} {\mathbb{N}}
\newcommand{\ZZ} {\mathbb{Z}}

\newcommand{\RR} {\mathbb{R}}
\newcommand{\CC} {\mathbb{C}}

\newcommand{\PP} {\mathbb{P}}
\renewcommand{\AA} {\mathbb{A}}

\newcommand {\shL}  {\mathcal{L}}

\newcommand {\shS}  {\mathcal{S}}
\newcommand {\shT}  {\mathcal{T}}

\newcommand {\shX}  {\mathcal{X}}
\newcommand {\shY}  {\mathcal{Y}}
\newcommand {\shZ}  {\mathcal{Z}}

\renewcommand {\Bar}  {\operatorname{Bar}}

\newcommand {\dual} {\vee}

\newcommand {\Hom}  {\operatorname{Hom}}

\newcommand {\Int}  {\operatorname{Int}}

\renewcommand{\O}  {\mathcal{O}}

\renewcommand{\P}  {\mathscr{P}}

\newcommand {\Spec} {\operatorname{Spec}}

\newcommand {\T} {\shT}
\newcommand {\X} {\shX}
\newcommand {\Y} {\shY}
\newcommand {\Z} {\shZ}

\newcommand {\momega} {\quad\hbox{\rput[r]{180}(-.3,.1){$\omega$}}}
\newcommand {\msig} {\quad\hbox{\rput[r]{180}(-.3,.1){$\sigma$}}}
\newcommand {\mtau} {\quad\hbox{\rput[r]{180}(-.3,.1){$\tau$}}}
\newcommand {\mbeta} {\quad\hbox{\rput[r]{180}(-.3,.1){$\beta$}}}
\newcommand {\R}  {\mathscr{R}}
\def\mydate{\ifcase\month \or January\or February\or March\or
April\or May\or June\or July\or August\or September\or October\or 
November\or December\fi \space\number\day,\space\number\year}

\begin{document}
\def\mapright#1{\smash{
  \mathop{\longrightarrow}\limits^{#1}}}
\def\mapleft#1{\smash{
  \mathop{\longleftarrow}\limits^{#1}}}
\def\exact#1#2#3{0\rightarrow#1\rightarrow#2\rightarrow#3\rightarrow0}
\def\mapup#1{\Big\uparrow
   \rlap{$\vcenter{\hbox{$\scriptstyle#1$}}$}}
\def\mapdown#1{\Big\downarrow
   \rlap{$\vcenter{\hbox{$\scriptstyle#1$}}$}}
\def\dual#1{{#1}^{\scriptscriptstyle \vee}}
\def\invlim{\mathop{\rm lim}\limits_{\longleftarrow}}

\input epsf.tex
\title[Toric Degenerations and Batyrev-Borisov Duality]{Toric Degenerations
and Batyrev-Borisov Duality}

\author{Mark Gross} \address{Department of Mathematics,
UCSD, 9500 Gilman Drive, La Jolla, CA 92093-0112, USA}
\email{mgross@math.ucsd.edu}
\thanks{This work was partially supported by NSF grant
0204326.}

\maketitle
\bigskip
\section*{Introduction.}

In \cite{AffI}, Bernd Siebert and I introduced the notion of a
{\it toric degeneration} of Calabi-Yau varieties. The initial goal is
to produce a method of constructing mirror pairs which combines the
Strominger-Yau-Zaslow (differential geometric) approach to mirror symmetry
and the older Batyrev-Borisov (algebro-geometric) approach to mirror
symmetry. Our belief is that in doing so we will produce a new, much more
powerful mechanism for constructing mirror pairs while gaining insight
into the underlying structures of mirror symmetry.

Sketching the approach here, we showed in \cite{AffI} 
how given a toric degeneration of Calabi-Yau varieties $\X\rightarrow\shS$,
one can associate to this degeneration a dual intersection complex $B$, which
is an integral affine manifold with singularities. Here $\shS$
is always a one-dimensional non-singular base. (See \S 1 of this paper
for a quick review of this construction). In addition, if the family
$\X\rightarrow\shS$  is polarized with a relatively ample divisor, then
$B$ comes along with a convex, piecewise linear multi-valued function
$\varphi$. We then defined a notion of discrete Legendre transform,
giving us a new affine manifold with singularities $\check B$ and a new
function $\check\varphi$. We can interpret $\check B$ as the dual intersection
complex of a new degeneration $\check\X\rightarrow\shS$. This latter part of the
program, at the time of this writing, is still not complete due to the
question of existence of smoothings, although we do not anticipate major
difficulties. In \cite{AffI}, we showed how to construct the degenerate fibre
$\check\X_0$, along with a log structure on $\check\X_0$ which is crucial
information for constructing a smoothing of $\check\X_0$. Once one knows
the log scheme $\check\X_0$ can be smoothed, then the smoothing can be
interpreted as an element of the mirror family.

Even without knowing whether the smoothing exists, this approach yields
important information about the SYZ approach. Indeed, $B$ contains
crucial information about the topology of the general fibre $\X_t$.
As $B$ is an integral affine manifold with singularities, there is an open
set $B_0\subseteq B$ which carries an integral affine structure, with
$\Gamma:=B\setminus B_0$ a locally finite union of locally closed submanifolds
of codimension at least two. Now we can define a topological space 
$X(B_0)$ as follows. Let $y_1,\ldots,y_n$ be local affine coordinates
on $B_0$, and define a local system $\Lambda\subseteq\T_{B_0}$ to be
the family of lattices generated by $\partial/\partial y_1,\ldots,
\partial/\partial y_n$. Because of the integrality of the affine structure,
this family of lattices is well-defined. We then define 
$X(B_0):=\T_{B_0}/\Lambda$; this is a torus bundle over $B_0$.
(This manifold carries a complex structure also, see \cite{Announce} for
details.) In addition, given an element ${\bf B}\in H^1(B_0,\Lambda\otimes(\RR/
\ZZ))$, one can twist the torus bundle by this element, obtaining
a new torus bundle $X(B_0,{\bf B})\rightarrow B_0$. In \cite{Torus}, we prove

\begin{theorem} 
\label{torusbundle}
Let $\X\rightarrow\shS$ be a toric degeneration, with
dual intersection complex $B$. Then
\item{(a)} For $t\not=0$, $\X_t$ is a topological compactification of
$X(B_0,{\bf B})$ for some ${\bf B}\in H^1(B_0,\Lambda\otimes (\RR/\ZZ))$.
\item{(b)} If $B$ is simple (see \cite{AffI}, Definition 1.60), 
then there exists a ${\bf B}\in H^1(B,i_*\Lambda\otimes(\RR/\ZZ))$
(with $i:B_0\hookrightarrow B$ the inclusion) and an algorithm
for compactifying $X(B_0,{\bf B})\rightarrow B_0$ to
a torus fibration $X(B,{\bf B})\rightarrow B$ so that $X(B,{\bf B})$
is homeomorphic to $\X_t$. 
\end{theorem}

If we are given a pair of toric degenerations $\X\rightarrow
\shS$ and $\check\X\rightarrow\shS$ whose dual intersection complexes are
related by a discrete Legendre transform, then the torus
bundles $X(B_0,{\bf B})\rightarrow B_0$ and $X(\check B_0,\check{\bf B})
\rightarrow \check B_0$ are in fact dual. Thus the above theorem
gives us a very general tool for proving topological forms of the SYZ
conjecture, so this approach is relevant even if one is not interested
in producing new mirror pairs, or not yet able to do so.

This paper should be regarded as an extended example of our construction,
though in fact it was part of the original motivation for the construction.
We consider mirror pairs constructed by the Batyrev-Borisov construction.
In this case, we already know of the existence of the mirror pairs, so we
are checking several things. First, we wish to know that our construction
gives the same answer as the Batyrev-Borisov construction. Second,
we wish to know that a topological form of the SYZ construction holds
for Batyrev-Borisov mirror pairs. This second fact follows from the first
by the above theorem, so this paper is concerned with considering the first
issue.

It is easy to construct examples of toric degenerations from the Batyrev-Borisov
construction. For brevity in the introduction, let us consider only
the toric hypersurface case. Let $M=\ZZ^n$, $M_{\RR}=M\otimes\RR$,
$\Delta\subseteq M_{\RR}$ be a reflexive polytope with $0$ as unique interior
integral point with reflexive dual $\Delta^*$. 
This defines a projective toric variety $(\PP_{\Delta},
\O_{\PP_{\Delta}}(1))$. 
The toric variety
$\PP_{\Delta}$ is defined by a fan $\Sigma$ in $N_{\RR}=\dual{M_{\RR}}$
consisting of cones over faces of $\Delta^*$. Similarly, $\Delta^*$ 
defines a projective toric variety $(\PP_{\Delta^*},\O_{\PP_{\Delta^*}}(1))$
defined by a fan $\check\Sigma$ in $M_{\RR}$.
If we take a general section $s$ of 
$\O_{\PP_{\Delta}}(1)$, and let $s_0$ be the section of $\O_{\PP_{\Delta}}(1)$
corresponding to $0\in\Delta$, then $ts+s_0=0$ defines a family $\X$ of
Calabi-Yau hypersurfaces in $\shS\times\PP_{\Delta}$, where $t$ is a coordinate
on $\shS$. For $t=0$, we obtain a degenerate hypersurface consisting of the
union of all toric divisors of $\PP_{\Delta}$. In particular 
$\X\rightarrow\shS$ is a toric degeneration. Now, it is very easy to
see, and was shown in \cite{AffI}, that the dual intersection complex
is just an affine structure (with singularities) 
on $\partial\Delta^*$, and that if we polarize the family with 
$\O_{\PP_{\Delta}}(1)$, then the discrete Legendre transform immediately
produces the dual intersection complex for a family of hypersurfaces
arising similarly 
from the dual $\Delta^*$. However, this is too weak a statement
for us, because these dual intersection complexes are probably never
simple, so we can't apply part (b) of Theorem \ref{torusbundle}.
Simplicity is essentially,
as we shall see, an ampleness statement (Theorem \ref{simplicity}). We need a more general
form of degeneration to get simplicity. So we use a standard trick,
already used in \cite{Zharkov} in the context of the SYZ conjecture, 
and probably first introduced by Viro. Choose a height function
$h:\Delta\cap M\rightarrow\ZZ$, and consider the degeneration
\begin{equation}
\label{degeneq}
s_0+\sum_{m\in\Delta\cap M} t^{h(m)}z^m=0.
\end{equation}
With suitable restriction on $h$,
when $t=0$ we again obtain just the equation $s_0=0$.
However, the total family $\X$ defined by this equation is very singular,
and toric techniques are required to resolve the singularities in order
to obtain a toric degeneration. 

More specifically, we will insist that $h$ be a restriction of an integral 
piecewise linear function on a subdivision of $\check\Sigma$ giving
a partial crepant resolution $\check\pi:\tilde\PP_{\Delta^*}\rightarrow
\PP_{\Delta^*}$ of $\PP_{\Delta^*}$. Such a piecewise linear function
determines a divisor $H_{\Delta^*}$ on $\tilde\PP_{\Delta^*}$, and we 
insist that $H_{\Delta^*}$ should be ample and $H_{\Delta^*}+
K_{\PP_{\Delta^*}}$ be nef. 

In this case, we show that the family defined by (\ref{degeneq}) is
birationally equivalent to a toric degeneration $\X\rightarrow\shS$,
with general fibre being a Calabi-Yau hypersurface in $\tilde\PP_{\Delta}$,
where $\pi:\tilde\PP_{\Delta}\rightarrow\PP_{\Delta}$ is any partial 
crepant projective resolution. From this toric degeneration we obtain
a dual intersection complex $B$. We show that if $\pi$ and $\check\pi$
are MPCP (maximal partial crepant projective) resolutions of
$\PP_{\Delta}$ and $\PP_{\Delta^*}$, then $B$ is simple. Of course, 
in any event,
we give an explicit combinatorial description of $B$ also.

The issue of duality is a bit more delicate. We show the following.
If $H_{\Delta}$ is an ample polarization on $\tilde\PP_{\Delta}$ induced
by a piecewise linear function $\check h$, (with
$H_{\Delta}+K_{\PP_{\Delta}}$ nef for symmetry), we obtain a polarization on
the generic fibre of $\X\rightarrow\shS$. However, we can not always show
that this polarization extends to $\X$. This seems to be a rather 
delicate combinatorial problem. However, if it does extend, then we
can perform a discrete Legendre transform on $B$, obtaining a dual $\check
B$, and we can show this is the dual intersection complex of a toric
degeneration $\check\X\rightarrow\shS$ obtained by interchanging $\Delta, h$
with $\Delta^*, \check h$. Thus the discrete Legendre transform provides
a more sophisticated form of Batyrev (-Borisov) duality. The only difficulty
is the existence of the polarization on $\X$. We obtain a weak
existence result, showing that there exists an $m_0$ such that
$m_0H_{\Delta}-n K_{\PP_{\Delta}}$ extends
to a polarization on $\X$ for $n$ sufficiently large. This is generally
sufficient for our purposes, for if $H_{\Delta}$ is ample on $\tilde\PP_{\Delta}$,
so is $m_0 H_{\Delta}-nK_{\PP_{\Delta}}$.

Of course, all of this is done much more generally in the complete intersection
case, which presents some new difficulties with the combinatorics; this
paper would be much shorter if we only treated the hypersurface case.

There has been important previous work covering the hypersurface
case. W.-D. Ruan in \cite{Ruan} gave a combinatorial description of the SYZ
fibration for hypersurfaces in toric varieties using his symplectic flow
argument; he obtains essentially the same description that we do, though he
doesn't use the language of affine manifolds. However, his proof of
duality is hard to understand. Haase and Zharkov give an
elegant explicit construction of the affine manifolds in the hypersurface
case, in which duality is self-evident. They have independently
generalized their construction to the complete intersection case \cite{HZ2}.

The structure of this paper is as follows. As we try to be as 
self-contained as possible, so that the reader need only consult
\cite{AffI} from time to time, in \S 1 we review the notion of toric 
degeneration and dual intersection complex. However,
notation of \cite{AffI}, \S 1 will be used throughout,
so some familiarity with that paper will be useful. 
In \S 2, we consider the basic
case of the Batyrev-Borisov construction, in which the polarizations
used are the anti-canonical ones. While covering this case first
may seem redundant, it is much simpler and much of the combinatorics
is needed in the general case, which is covered in \S 3. In \S 4,
we consider a worked example, one which gives a particularly nice
example of a Strominger-Yau-Zaslow fibration which was recently studied by Kovalev.
In this case, we see how the language of toric degenerations gives
a much more flexible picture than the Batyrev-Borisov construction
alone does.

I would like to thank Bernd Siebert, with whom the overall project of
studying mirror symmetry via toric degenerations has been joint work,
and Alexei Kovalev and Ilia Zharkov for useful conversations.

\section{Review}

We recall the definition of toric degeneration from \cite{AffI}.

\begin{definition}
\label{toric degen}
Let $R$ be a discrete valuation ring and $k$-algebra with algebraically closed
residue class field $k$. A toric degeneration of Calabi-Yau
varieties over $R$ is a proper normal algebraic space $\X$ flat over
$\shS:=\Spec R$ satisfying the following properties:
\item{(1)} The generic fibre $\X_{\eta}$ is an irreducible normal
variety over $\eta$.
\item{(2)} $\X_0$ is reduced, and 
if $\nu:\tilde\X_0\to\X_0$ is the normalization,
then $\tilde\X_0$ is a disjoint union of toric varieties,
the conductor scheme $C\subseteq\tilde\X_0$ is reduced
and the map $C\to\nu(C)$ is unramified and generically
two-to-one. The square
\[\begin{CD}
C@>>> \tilde\X_0\\
@VVV @VV{\nu}V\\
\nu(C)@>>> \X_0
\end{CD}\]
is cartesian and cocartesian.
\item{(3)} $\X$ (and hence $\X_0$) is Gorenstein, and
the conductor locus $C$ restricted to each irreducible component
of $\tilde\X_0$ is the union of all toric Weil divisors.
\item{(4)} There exists a closed subset $Z\subseteq\X$ of relative
codimension $\ge 2$ such that $Z$ satisfies the following properties:
$Z$ does not contain the image under $\nu$ of any toric stratum of
$\tilde\X_0$, and for any geometric point $\bar x\to \X\setminus Z$,
there is an \' etale neighbourhood $U_{\bar x}\to \X\setminus Z$ of
$\bar x$, an affine toric variety $Y_{\bar x}$, a regular function
$f_{\bar x}$ on $Y_{\bar x}$ given by a monomial, a choice of
uniformizing parameter of $R$ giving a map $k[\NN]\rightarrow R$, and
a commutative diagram
\[
\begin{matrix}
U_{\bar x}&\mapright{}&Y_{\bar x}\\
\mapdown{f|_{U_{\bar x}}}&&\mapdown{f_{\bar x}}\\
\Spec R&\mapright{}&\Spec k[\NN]
\end{matrix}
\]
such that the induced map $U_{\bar x}\to \Spec R\times_{\Spec k[\NN]}
Y_{\bar x}$ is smooth. Furthermore, $f_{\bar x}$ vanishes on each
toric divisor of $Y_{\bar x}$.
\end{definition}

The slogan to remember here is that this definition roughly implies
that $\X_0$ is a union of toric varieties meeting along toric strata,
$\X_0$ has numerically trivial canonical class, and away from $Z$,
$\X\rightarrow\shS$ is described \'etale locally as a monomial
map on an affine toric variety.

We recall how to construct the dual intersection complex of a toric
degeneration $\X\rightarrow\shS$. We will make the simplifying assumption
in this review that all irreducible components of $\X_0$ are normal.
We do this here as this will be the case with all examples constructed
in this paper.

Write
$\tilde\X_0=\coprod X_i$, with $\nu:\tilde \X_0\to \X_0$ the normalization.
The set of \emph{strata} of $\X_0$ is the set
\[
Strata(\X_0):=\{\nu(S)|\hbox{$S$ is a toric stratum of $X_i$ for some $i$}\}.
\]

Associate to each strata of $\X_0$ a lattice polytope $P(S)$
as follows. If $S\in Strata(\X_0)$, with generic point $\eta$, 
then by Condition (4) of Definition \ref{toric degen}, there is
a toric variety $Y_{\bar\eta}$ defined by a convex rational polyhedral
cone $\sigma\subseteq M_{\RR}$, where $M$ is a lattice
and $M_{\RR}=M\otimes\RR$,
so that $Y_{\bar\eta}=\Spec k[\dual{\sigma}\cap N]$ with $N=\Hom(M,\ZZ)$.
In addition, there is an element $\rho\in\dual{\sigma}\cap N$
inducing the map $Y_{\bar\eta}\rightarrow \Spec k[\NN]$. Define
\[
P(S):=\{m\in\sigma|\langle\rho,m\rangle=1\}.
\]
One can show (see \cite{AffI}, \S 4) that $P(S)$ is a lattice polytope
of dimension equal to the codimension of $S$ in $\X_0$.

Furthermore, if $S_1\subseteq S_2$ are two strata with corresponding cones
$\sigma_1$ and $\sigma_2$, then toric geometry tells us that the faces
of $\sigma_1$ are in one-to-one inclusion reversing correspondence
with toric strata containing $S_1$. In particular, $\sigma_2$ is
a face of $\sigma_1$ and $P(S_2)$ is naturally a face of $P(S_1)$.
We can now construct the dual intersection complex $B$ of $\X$ as a
union of polyhedra. Explicitly, if $S_1$ and $S_2$ are two zero-dimensional
strata with $S_3$ the minimal stratum containing $S_1$ and $S_2$, then we 
identify $P(S_1)$ and $P(S_2)$ along the faces $P(S_3)\subseteq P(S_1)$ and
$P(S_3)\subseteq P(S_2)$. In \cite{AffI}, Proposition 4.10, it was proved that
$B$ is an $n$-dimensional manifold, where $n=\dim\X_0$. 
Furthermore, $B$ comes along
with a polyhedral decomposition
\[
\P=\{P(S)|S\in Strata(\X_0)\}.
\]
Finally, $B$ can be given the structure of an integral affine manifold with 
singularities by specifying a {\it fan structure} at each vertex of
$\P$. For each $v$ a vertex of $\P$, there is a corresponding irreducible
component $X_v$ of $\X_0$, defined by a fan $\Sigma_v\subseteq\RR^n$. 
Now each cone of $\Sigma_v$ corresponds to a toric stratum
of $X_v$, and hence to an element of $\P$. This gives a 
combinatorial correspondence between cones of $\Sigma_v$ and
cells of $\P$ containing $v$. This correspondence is a bijection in the
case that all irreducible components of $\X_0$ are normal. Let $W_v$
be the union of interiors of all simplices in the first barycentric
subdivision of $\P$ containing $v$ (the star of $v$). Then it is not
difficult to see that one can define a map $\psi_v:W_v\rightarrow
\RR^n$ which is a homeomorphism onto its image, takes each cell
of $\P$ containing $v$ into the corresponding cone of $\Sigma_v$,
and is integral affine on the interiors of each maximal cell of $\P$
containing $v$. Together with the obvious charts
\[
\psi_S:Int(P(S))\hookrightarrow \AA_S:=\{m\in M_{\RR}|\langle \rho,m\rangle=1\}
\]
for $S$ a minimal stratum, we obtain an integral affine structure on
$B_0:=B\setminus\Gamma$, where $\Gamma$ is the union of all simplices of 
$\Bar(\P)$, the first barycentric subdivision of $\P$,
not containing a vertex of $\P$ or intersecting the
interior of a maximal cell of $\P$. This discriminant locus is usually
much larger than necessary.

Recall that as $B_0$ carries an integral affine structure, we can define
a local system $\Lambda\subseteq\T_{B_0}$ which is the family of
integral lattices locally generated by $\partial/\partial y_1,\ldots,
\partial/\partial y_n$, where $y_1,\ldots,y_n$ are local integral
affine coordinates. In \cite{AffI}, Proposition 1.27, 
it is shown that if the monodromy of
$\Lambda$ about a simplex of $\Gamma$ is trivial, the affine structure
of $B_0$ can be extended across this simplex. We shall use this latter
fact in \S 4 to obtain a minimal discriminant locus in the example
given there.
\bigskip

We will not review the construction of the discrete Legendre transform
here, but rather refer the reader to either \cite{Announce}, \S 4 or
\cite{AffI}, \S 1.4. Briefly, given an integral affine manifold with
singularities $B$ with a polyhedral decomposition $\P$ and a multi-valued
strictly convex piecewise linear function $\varphi$ on $B$, the discrete
Legendre transform of the triple $(B,\P,\varphi)$ is a new triple
$(\check B,\check\P,\check\varphi)$ with $B=\check B$ topologically, but
with different affine structures. The affine structures are dual
in the sense that $\Lambda_{\check B}=\Hom(\Lambda_B,\ZZ)$. Thus the torus
bundles $X(B_0)\rightarrow B_0$ and $X(\check B_0)\rightarrow B_0$
are dual in the Strominger-Yau-Zaslow sense. 

If $\X\rightarrow\shS$ is a toric degeneration and $\shL$ is a relatively
ample line bundle on $\X$, we say we have a \emph{polarized toric degeneration}.
Restricting $\shL$ to each irreducible component of $\X_0$ gives an ample
line bundle on a toric variety, hence a strictly convex piecewise linear
function on the corresponding fan. This gives rise to a multi-valued piecewise
linear function $\varphi$ on the dual intersection complex $B$ of $\X
\rightarrow\shS$.

The discrete Legendre transform of $(B,\P,\varphi)$ in this case
can be seen as the \emph{intersection complex} of $\X\rightarrow\shS$.
This is obtained as follows. For every $k$-dimensional 
toric strata $S$ of $\X_0$,
$\shL|_S$, being an ample line bundle on $S$, has a Newton polytope,
a $k$-dimensional lattice polytope $\check P(S)\subseteq\RR^k$, well
defined up to integral affine transformation. If $S_1\subseteq S_2$, then we can
canonically identify $\check P(S_1)$ with a face of $\check P(S_2)$.
Making these identifications we obtain $\check B,\check\P$. This is
intuitively just a ``real'' drawing of the central fibre $\X_0$, obtained
by dividing each component of $\X_0$ by its canonical real torus action.

To obtain an affine sturcture on $\check B$, we proceed as in the construction
for $B$: each maximal $\check P(S)$ has a canonical affine structure
in its interior; for each zero-dimensional stratum $S$ we specify the
fan structure at the vertex $\check P(S)$ to be the normal fan to
$P(S)$. This yields $\check B$ as an affine manifold with singularities.

\section{The Batyrev-Borisov construction: the basic case}

We first recall the construction of Batyrev and Borisov. We will then
demonstrate that the Batyrev-Borisov construction is a case of the
discrete Legendre transform and gives rise to toric degenerations. This
will be done first with the anti-canonical polarization, and we will
cover the general polarization in the next section.

\begin{definition}
Let $M=\ZZ^n$, $N=\Hom(M,\ZZ)$,
and let $\Delta\subseteq M_{\RR}$ be a reflexive
polytope. 
Let $\Delta=\Delta_1+\cdots+\Delta_r$, $\Delta_i\subseteq
M_{\RR}$ be a Minkowski decomposition of $\Delta$ into lattice polytopes.
Let $\Sigma\subseteq N_{\RR}$
be the normal fan to $\Delta$
defining $\PP_{\Delta}$, and $\varphi:N_{\RR}\rightarrow
{\RR}$ be the integral (upper) convex $\Sigma$-piecewise linear function
corresponding to the anti-canonical sheaf. 
If $\Sigma(1)=\{e_1,\ldots,e_m\}$ (the set of integral generators of
one-dimensional cones of $\Sigma$),
$\varphi(e_i)=1$. Then $\Delta=\Delta_1+\cdots+\Delta_r$ yields a 
decomposition $\varphi=\varphi_1+\cdots+\varphi_r$, where $\varphi_i$
defines the divisor corresponding to $\Delta_i$. Here
$$\varphi_i(e_j)=-\inf \{\langle x,e_j\rangle|x\in\Delta_i\},$$
and
$$\Delta_i=\{x\in M_{\RR}|\langle x,y\rangle \ge -\varphi_i(y)\quad
\forall y\in N_{\RR}\}.$$
We say $\Delta=\Delta_1+\cdots+\Delta_r$ is a {\it nef-partition}
if $\varphi_i(e_j)
\in \{0,1\}$ for all $i,j$. 
\end{definition}

A nef-partition $\Delta=\Delta_1+\cdots+\Delta_r$ yields a decomposition
$\O_{\PP_{\Delta}}(1)\cong \shL_1\otimes\cdots\otimes\shL_r$, with
$\shL_i$ a semi-ample line bundle corresponding to $\Delta_i$. Recall
(essentially \cite{BB}, Proposition 4.15)

\begin{proposition}
\label{general}
Let $\Delta=\Delta_1+\cdots+\Delta_r$ be a
nef partition, and $s$ be a general section of
$\shL_1\oplus\cdots\oplus\shL_r$. Let $V\subseteq\PP_{\Delta}$ be the zero-locus
of $s$. If $V$ is non-empty and irreducible, then $V$ is
an $n-r$ dimensional variety, $\omega_V\cong\O_V$, and $V$ has only canonical
singularities.
\end{proposition}

\begin{definition}
Given a nef partition, let $\nabla_i$ be the convex hull of $0\in N_{\RR}$
and all $e_j$ such that $\varphi_i(e_j)=1$. Set $\nabla=\nabla_1+\cdots+
\nabla_r$. 
\end{definition}

\begin{proposition} If $\nabla^*,\Delta^*$
are the duals of $\nabla$ and $\Delta$ respectively, then
\begin{eqnarray*}
\nabla^*&=&Conv\{\Delta_1,\ldots,\Delta_r\}\\
\Delta^*&=&Conv\{\nabla_1,\ldots,\nabla_r\}.
\end{eqnarray*}
\end{proposition}

Proof. \cite{BB}, Theorem 4.10. \qed

We also will use

\begin{lemma}
\label{borlem}
$$\langle \Delta_j,\nabla_i\rangle \ge-\delta_{ji}.$$
\end{lemma}

Proof. \cite{Borisov}, Cor. 2.12. \qed

Thus, to summarize: 
We are initially given data 
\begin{eqnarray*}
\Delta&=&\Delta_1+\cdots+\Delta_r\\
\nabla&=&\nabla_1+\cdots+\nabla_r\\
\Delta^*&=&Conv\{\nabla_1,\ldots,\nabla_r\}\\
\nabla^*&=&Conv\{\Delta_1,\ldots,\Delta_r\}\\
\varphi&=&\varphi_1+\cdots+\varphi_r\\
\check\varphi&=&\check\varphi_1+\cdots+\check\varphi_r
\end{eqnarray*}
where $\varphi,\varphi_i$ are piecewise linear functions on the
normal fan $\Sigma$ to $\Delta$, $\varphi$ representing $-K_{\PP_{\Delta}}$,
and $\check\varphi$, $\check\varphi_i$ are piecewise linear functions
on the normal fan $\check\Sigma$ to $\nabla$, $\check\varphi$ representing
$-K_{\PP_{\nabla}}$. We have
\begin{eqnarray*}
\nabla_i&=&Conv(\{n\in\Sigma(1)|\varphi_i(n)=1\}\cup\{0\})\\
\Delta_i&=&Conv(\{m\in\check\Sigma(1)|\check\varphi_i(m)=1\}\cup\{0\}).
\end{eqnarray*}

We first work out the combinatorics of this duality from our point of
view. In the hypersurface case, there is a straightforward dual
relationship between faces of $\Delta$ and $\nabla=\Delta^*$. In the
complete intersection case this is more delicate.

\begin{definition} 
Let $\sigma^*\subset\Delta^*$ be a proper face. Let
$$\sigma=\{m\in\Delta|\hbox{$\langle m,n\rangle =-1$ for all $n\in
\sigma^*$}\}$$
be the dual face of $\Delta$. Let 
$$\sigma_i^*=\{n\in\sigma^*|\varphi_i(n)=1\}\subseteq\nabla_i$$
and set
$\msig=\sigma_1^*+\cdots+\sigma_r^*\subseteq\nabla$.
(We take $\msig$ to be empty if any $\sigma_i^*$ is empty). Let
$$\msig^*=\{m\in\nabla^*|\hbox{$\langle m,n\rangle =-1$ for all $n\in
\msig$}\},$$
and let 
$$\msig_i^*=\{m\in\msig^*|\check\varphi_i(m)=1\}\subseteq\Delta_i.$$
\end{definition}

\begin{lemma}
\label{corrbetweenfaces}
\item{(a)}
Let $\sigma^*\subset\Delta^*$ be a proper face, and suppose $\msig$ is
non-empty. Then
\begin{eqnarray*}
\sigma^*&=&Conv(\sigma_1^*,\ldots,\sigma_r^*)\\
\msig^*&=&Conv(\msig_1^*,\ldots,\msig_r^*)\\
\sigma&=&\msig_1^*+\cdots+\msig_r^*
\end{eqnarray*}
and $\msig$ is a proper face of $\nabla$.
\item{(b)} 
\begin{eqnarray*}
\msig_i^*&=&\{m\in\Delta_i|\hbox{$\langle m,n\rangle=-\varphi_i(n)$ 
for all $n\in\sigma^*$}\}\\
\sigma_i^*&=&\{n\in\nabla_i|\hbox{$\langle m,n\rangle=-\check\varphi_i(m)$ 
for all $m\in\msig^*$}\}
\end{eqnarray*}
\end{lemma}

Proof. (a) Clearly $Conv(\sigma_1^*,\ldots,\sigma_r^*)\subseteq
\sigma^*$. Conversely, let $n\in\sigma^*$ be a vertex. Then $\varphi_i(n)
=1$ for some $i$, so $n\in\sigma_i^*$. Thus $\sigma^*\subseteq 
Conv(\sigma_1^*,\ldots,\sigma_r^*)$. Similarly, 
\[
\msig^*=
Conv(\msig_1^*,\ldots,\msig_r^*).
\]

Next let $m\in\sigma$ be a vertex of $\sigma$. Since $\Delta=\Delta_1+
\cdots+\Delta_r$, we can write $m=m_1+\cdots+m_r$ where $m_i$
is a vertex of $\Delta_i$. Suppose $m_i=0$. If
$n\in\sigma^*_i\subseteq\sigma^*$ with $\varphi_i(n)=1$ (such exists by
the assumption that $\msig$ is non-empty), then 
$-1=\langle m,n\rangle =\sum_{j=1\atop i\not=j}^r \langle m_j,n\rangle
\ge 0$ by Lemma \ref{borlem}. Thus $m_i\not=0$, so each $m_i$
is a non-zero vertex of $\Delta_i$, and
$\check\varphi_i(m_i)=1$. Next, let
$n\in\msig$, $n=n_1+\cdots+n_r$ with $n_i\in\sigma_i^*\subseteq\sigma^*$.
Then
$$-1=\langle m,n_j\rangle=\sum_{i=1}^r \langle m_i,n_j\rangle
\ge -\sum_{i=1}^r\delta_{ij},$$
so in fact $\langle m_i,n_j\rangle =-\delta_{ij}$. Thus
$\langle m_i,n\rangle =-1$, so $m_i\in\msig^*$, and hence $m_i\in\msig^*_i$.
We conclude $\sigma\subseteq \msig_1^*+\ldots+\msig_r^*$.
(Note in particular, we have shown $\msig^*$ is non-empty if
$\sigma$ is non-empty.)

Conversely, let $m=m_1+\cdots+m_r\in\msig_1^*+\cdots+\msig_r^*$.
Then $m\in\Delta$. Let $n\in\sigma^*$ be a vertex, so $n\in \sigma_i^*$
for some $i$. Choose $n_j\in\sigma_j^*$ for $j\not=i$, and set $n_i=n$.
Then as $m_j\in\msig^*$, $-1=\langle m_j,\sum_{k=1}^r n_k\rangle
\ge -\sum_{i=1}^r \delta_{jk}$, so $\langle m_j,n_k\rangle =-\delta_{jk}$.
Thus $\langle m,n\rangle =\langle m,n_i\rangle =-1$. Since this holds
for every vertex $n$ of $\sigma^*$, $m\in\sigma$. Thus
$\sigma=\msig_1^*+\cdots+\msig_r^*$.

Finally, as noted above, $\msig^*$ is non-empty, and hence $\msig$ is
contained in a proper face of $\nabla$. Let $\msig'$ denote
the minimal face of $\nabla$ containing $\msig$, i.e. the face dual
to $\msig^*$. Then interchanging $\Delta$ with $\nabla$ and 
$\nabla^*$ with $\Delta^*$, the statement $\sigma=\msig_1^*+\cdots+
\msig_r^*$ proved above becomes $\msig'=\sigma_1^*+\cdots+\sigma_r^*$.
Thus $\msig=\msig'$.

(b) We just prove the first statement. If $m=m_1+\cdots+m_r$ with
$m_i\in\msig_i^*$, $n_j\in\sigma_j^*\subseteq\sigma^*$,
then we saw that $\langle m_i,n_j\rangle=-\delta_{ij}=-\varphi_i(n_j)$,
from which it follows that $\msig_i^*\subseteq\{m\in\Delta_i|
\langle m,n\rangle=-\varphi_i(n)\quad\forall n\in\sigma^*\}$. 
Conversely, let $m\in\Delta_i$ be a vertex with $\langle m,n\rangle
=-\varphi_i(n)$ for all $n\in\sigma^*$. Since $\msig\not=\emptyset$,
this means $m\not=0$, and $\check\varphi_i(m)=1$. Thus if 
$m\in\msig^*$, $m\in\msig^*_i$ as desired. Now if $n_j\in\sigma_j^*
\subseteq\sigma^*$, $n=n_1+\cdots+n_r\in\msig$, and $\langle
m,n\rangle=-\sum_j \varphi_i(n_j)=-1$, so $\langle m,n\rangle
=-1$ for all $n\in\msig$, so $m\in\msig^*$.
\qed

It is clear from this that given any one of $\sigma$, $\sigma^*$,
$\msig$ or $\msig^*$, we obtain the other three. We now state
the following combinatorial duality, which follows immediately
from the Lemma.

\begin{corollary}
\label{correspondence}
Let 
$$\R_{\Delta^*}=\{\sigma^*\subset\Delta^*| \msig\not=\emptyset\}$$
and
$$\R_{\nabla^*}=\{\msig^*\subset\nabla^*| \sigma\not=\emptyset\}$$
be subsets of the set of faces of $\Delta^*$ and $\nabla^*$
respectively. Let
$$\P_{\nabla}=\{\msig\subset\nabla|\sigma^*\in\R_{\Delta^*}\}$$
and
$$\P_{\Delta}=\{\sigma\subset\Delta|\msig^*\in\R_{\nabla^*}\}.$$
Then the maps $\R_{\Delta^*}\rightarrow\P_{\nabla}$ and $\R_{\nabla*}
\rightarrow\P_{\Delta}$ given by $\sigma^*\mapsto\msig$ and $\msig^*
\mapsto \sigma$ are order preserving bijections. The maps
$\R_{\Delta^*}\rightarrow \R_{\nabla^*}$ given by $\sigma^*\mapsto
\msig^*$ and $\P_{\nabla}\rightarrow\P_{\Delta}$ given by
$\sigma\mapsto\msig$ are order reversing bijections.
\end{corollary}

\begin{proposition}
\label{dims}
If $\msig\in\P_{\nabla}$, any face of $\msig$ is also in $\P_{\nabla}$,
and dually, for $\sigma\in\P_{\Delta}$. Also, if $\msig\in\P_{\nabla}$, then
\begin{eqnarray*}
\dim\msig&=&\dim\sigma^*-r+1\\
\dim\sigma&=&\dim\msig^*-r+1\\
\dim\msig&=&(\dim M_{\RR}-r)-\dim\sigma
\end{eqnarray*}
\end{proposition}

Proof. Let $\mtau$ be a face of $\msig$. Then we obtain $\mtau^*\supseteq
\msig^*$, and hence $\mtau_i^*\supseteq \msig_i^*$, and hence
$\mtau_i^*$ is non-empty for each $i$. Taking $\tau=\mtau_1^*+\cdots
+\mtau_r^*\supseteq\sigma$, we obtain $\tau^*\subseteq\sigma^*$, and by
the Lemma, $\mtau=\tau_1^*+\cdots+\tau_r^*$. Thus $\mtau\in\P_{\nabla}$.

It is enough to prove the first dimension statement. The second will be
dually true and then $\dim\msig=\dim\sigma^*-r+1
=(\dim M_{\RR}-\dim\sigma-1)-r+1$. To show the first dimension statement,
we proceed as follows. First, if $\dim\msig=0$, then $\dim\sigma_i^*=0$
for each $i$, and $\sigma^*$ is a convex hull of vertices $v_1,\ldots,v_r$,
with $\sigma_i^*=\{v_i\}$. Furthermore, since $\varphi_i(v_j)=\delta_{ij}$,
$v_1,\ldots,v_r$ are linearly independent, so $\dim\sigma^*=r-1$.

Now let $\msig\in\P_{\nabla}$ be arbitrary, with $\dim\msig=p$. Choose
a vertex $\mtau_0$ of $\msig$, and a vertex $\tau_q$ of $\sigma$,
(here $q$ will be an index we will set in a moment equal to $\dim M_{\RR}-r$). 
Then $\tau_q^*$ is a maximal face of $\Delta^*$ containing $\sigma^*$,
and hence $\mtau_q$ is a maximal element of $\P_{\nabla}$ containing
$\msig$. On the other hand, arguing dually as above, since $\tau_q$ is
a vertex, $\dim \mtau_q^*=r-1$, so $\dim \mtau_q=\dim M_{\RR}-r$.
Now take a strictly ascending 
chain of faces $\mtau_0\subset\cdots\subset\mtau_q$ with $\mtau_p=\msig$,
so that $\dim\mtau_i=i$. We then obtain a strictly ascending chain of
faces $\tau_0^*\subset\cdots\subset\tau_q^*$ with $\tau_p^*=\sigma^*$.
Since $\dim\tau_0^*=r-1$ and $\dim \tau_q^*=\dim M_{\RR}-1$, it follows
that $\dim\tau_i^*=r-1+i$ and $\dim\sigma^*=r-1+\dim\msig$, as desired.
\qed

\begin{definition} 
Let
\begin{eqnarray*}
B_{\nabla}&=&\bigcup_{\msig\in\P_{\nabla}} \msig\subseteq\nabla\\
B_{\Delta}&=&\bigcup_{\sigma\in\P_{\Delta}} \sigma\subseteq\Delta
\end{eqnarray*}
so that $\P_{\nabla}$ and $\P_{\Delta}$ give polyhedral decompositions
of $B_{\nabla}$ and $B_{\Delta}$ respectively, with first barycentric
subdivisions $\Bar(\P_{\nabla})$ and $\Bar(\P_{\Delta})$.
We define a discriminant locus of $B_{\nabla}$, $\Gamma_{\nabla}$,
as the union
of all simplices of $\Bar(\P_{\nabla})$ not containing a vertex of $\P_{\nabla}$
or intersecting the interior of a maximal cell of $\P_{\nabla}$. We
do the same dually to define $\Gamma_{\Delta}\subseteq B_{\Delta}$.

We define an affine structure on $B_{\nabla}\setminus\Gamma_{\nabla}$ as 
follows.
For every maximal cell $\msig$ of $B_{\nabla}$, let ${\AA}_{\msig}$ be the
affine subspace of $N_{\RR}$ spanned by $\msig$. Then we define
an affine chart on $\Int(\msig)$ by  $\psi_{\msig}:\Int(\msig)\hookrightarrow
{\AA}_{\msig}$. Secondly, for $v\in\P_{\nabla}$ a vertex, we have 
$$v=\sum_{i=1}^r v_i$$
for some $v_i\in \Sigma(1)$. 
Let $W_v$ be the
(open) star of $v$ in $\Bar(\P_{\nabla})$, i.e. the union of interiors of all
simplices of $\Bar(\P_{\nabla})$ containing $v$.
We then define an affine chart
$$\psi_v:W_v\rightarrow N_{\RR}/Span(v_{1},\ldots,v_{r})$$
simply via projection.

Similarly we define an affine structure on $B_{\Delta}\setminus\Gamma_{\Delta}$.
\qed
\end{definition}

If $\sigma\subseteq \RR^n$, we always write 
\[
C(\sigma)=\{rm| r\ge 0, m\in \sigma\}
\]
for the cone over $\sigma$.

\begin{lemma} 
\label{fan combo}
\item{(a)} Let $\sigma^*\in\R_{\Delta^*}$. Let $m_1,\ldots,m_r
\in M_{\RR}/C(\sigma^*)^{\perp}$ be such that $\langle m_i,n\rangle
=\varphi_i(n)$ for all $n\in\sigma^*$, and write 
$$V_{\sigma}=Span\{m_i-m_j| 1\le i,j \le r\}\subseteq M_{\RR}/C(\sigma^*)
^{\perp}.$$
This gives a subspace $V_{\sigma}^{\perp}\subseteq\RR C(\sigma^*)\subseteq
N_{\RR}$. Then 
$$C(\msig)=V_{\sigma}^{\perp}\cap C(\sigma^*).$$
\item{(b)} If $\tau^*\subseteq\sigma^*$, $\tau^*,\sigma^*\in \R_{\Delta^*}$,
then there is an isomorphism $\RR C(\msig)/\RR C(\mtau)\cong
\RR C(\sigma^*)/\RR C(\tau^*)$ inducing an isomorphism of cones
$$(C(\msig)+\RR C(\mtau))/\RR C(\mtau)
\cong
(C(\sigma^*)+\RR C(\tau^*))/\RR C(\tau^*).$$
\end{lemma}

Proof. (a) Since $\langle m_i,\sigma_j^*\rangle =\delta_{ij}$, $m_i=m_j$
on $\msig=\sigma_1^*+\cdots+\sigma_r^*$. Thus $\msig\subseteq
V_{\sigma}^{\perp}\cap C(\sigma^*)$. 
Conversely, if $n\in V_{\sigma}^{\perp}\cap C(\sigma^*)$,
then either $\langle m_i,n\rangle=0$ for all $i$, in which case $n=0$,
or we can multiply $n$ by a positive real number so that we can
assume $\langle m_i,n\rangle=1$ for all $i$. We wish to prove that
after doing so, $n\in \msig$. For this, it is enough to show
$\langle m,n\rangle =-1$ if $m\in\msig^*$ and $\langle m,n\rangle\ge -1$
if $m\in\nabla^*$. But $\msig^*=Conv(\msig_1^*,\ldots,\msig_r^*)$, and
if $m\in \msig_i^*$, then $m$ defines the function $-\varphi_i$ on 
$C(\sigma^*)$ by Lemma \ref{corrbetweenfaces}, 
(b), and so $\langle m,n\rangle =-1$.
Thus if $m\in\msig^*$, $\langle m,n\rangle =-1$. If $m\in\Delta_i$,
then $\langle m,n\rangle\ge -\varphi_i(n)=-1$, so if $m\in\nabla^*
=Conv(\Delta_1,\ldots,\Delta_r)$, $\langle m,n\rangle\ge -1$. Thus
$n\in\msig$ as desired.

(b) Note that $\RR C(\msig)=V_{\sigma}^{\perp}$: certainly by (a)
$\RR C(\msig)\subseteq V_{\sigma}^{\perp}$, and $\dim C(\msig)
=\dim C(\sigma^*)-r+1$ by Proposition \ref{dims} while 
$\dim V_{\sigma}^{\perp}=\dim \RR C(\sigma^*)-(r-1)=\dim C(\msig)$. Thus
the two vector spaces coincide. Now we have a map $V_{\sigma}^{\perp}/
V_{\tau}^{\perp}\rightarrow \RR C(\sigma^*)/\RR C(\tau^*)$, which is
injective since $V_{\tau}^{\perp}=V_{\sigma}^{\perp}\cap \RR C(\tau^*)$,
and is surjective,
as can be easily seen by calculating dimensions. Since $C(\msig)\subseteq
C(\sigma^*)$, clearly we then get an inclusion of cones
$(C(\msig)+\RR C(\mtau))/\RR C(\mtau)
\subseteq
(C(\sigma^*)+\RR C(\tau^*))/\RR C(\tau^*)$. On the other hand, if $n\in 
C(\sigma^*)$, there exists $n'\in \RR C(\tau^*)$ such that 
$n+n'\in \RR C(\msig)$. Now, if there is an $n''\in C(\mtau)$
such that $n''\not\in\partial C(\tau^*)$, then for some sufficiently
large $r\in\RR$, $n'+rn''\in C(\tau^*)$. However, it is easy to see
that $C(\mtau)\not\subseteq \partial C(\tau^*)$, hence we have
$n+n'+rn''\in C(\msig)$. Thus the inclusion of cones is 
a surjection. \qed

\begin{proposition} The charts $\psi_v$ and $\psi_{\sigma}$
define an integral affine structure on $B_{\nabla}\setminus\Gamma_{\nabla}$,
making $B_{\nabla}$ an integral affine manifold with singularities
of dimension $\dim M_{\RR}-r$.
\end{proposition}

Proof. The open sets of the form $\Int(\sigma)$ and $W_v$ clearly
cover $B_{\nabla}\setminus\Gamma_{\nabla}$, so to show we obtain an
integral affine structure, we need to show $\psi_v$ is a $C^0$ immersion
and that the transition maps are integral affine. Let $v=\sum_{j=1}^r v_j$
be a vertex of $\P_{\nabla}$ 
and let $\msig$ be any maximal cell containing
$v$ and ${\AA}_{\msig}$ the affine space spanned by it. 
The functions $\varphi_i$ are linear on the cone over $\sigma^*$, and hence are
induced by linearly independent elements $m_i\in M$ on $\sigma^*$.
The space
${\AA}_{\msig}$ is then given by
$${\AA}_{\msig}=\{x\in N_{\RR}|
\hbox{$\langle m_i,x\rangle=1$ for $i=1,\ldots,r$}\},$$
as follows from Lemma \ref{fan combo}.

Also
$$\langle m_i,\sum_{j=1}^r a_jv_j\rangle =a_i,$$
so if $\sum_{j=1}^r a_jv_j\in {\AA}_{\msig}$, we have
$a_j=1$, for all $j$. Thus
$${\AA}_{\msig}\cap Span(v_1,\ldots,v_r)=\{v\}.$$
In particular, the projection
$${\AA}_{\msig}\rightarrow N_{\RR}/Span(v_1,\ldots,v_r)$$
is an affine isomorphism preserving integral structures, i.e.
mapping ${\AA}_{\msig}\cap N$ isomorphically to 
$N/(Span(v_{1},\ldots,v_{r})\cap N)$. Thus the transition
maps are integral affine.

We then have to show that $\psi_v$ is 1-1. First let's look at the tangent
wedge to $\psi_v(\msig)$ at $0=\psi_v(v)$. Let $\tau^*=Conv(v_1,\ldots,v_r)$.
Then the tangent cone to 
$\psi_v(\msig)$ is precisely the cone $(C(\msig)+\RR v)/\RR v$, equal
to $(C(\sigma^*)+\RR C(\tau^*))/\RR C(\tau^*)$
of the fan $\Sigma(C(\tau^*))$ by Lemma \ref{fan combo}, (b),
(see \cite{AffI}, Definition 1.37 for notation). Since the maximal facets
of $\P_{\nabla}$ containing $v$ are in one-to-one correspondence with maximal
cones of $\Sigma$ containing $C(\tau^*)$, we see that the tangent wedges to
all the $\psi_v(\msig)$'s fit together to form the fan $\Sigma(C(\tau^*))$.
Thus in particular $\psi_v$ is a ($C^0$) immersion.
\qed
\bigskip

The computation of monodromy can be done via a calculation of the sort
carried out by Ruan in \cite{Ruan} or Haase and Zharkov in \cite{HZ}.

\begin{proposition} 
\label{monodromy1}
Let $v$ and $v'$ be two vertices of $\P_{\nabla}$
and let $\msig$ and $\msig'$ be two
maximal faces of $\P$ containing $v$ and $v'$. Let $\gamma$ be
a simple loop based at $v$, passing successively into $\Int(\msig)$,
through $v'$, into $\Int(\msig')$ and back to
$v$. Write $v=v_1+\cdots+v_r$ and $v'=v_1'+\cdots+v_r'$,
for $v_j,v_j'\in\Sigma(1)$, and let $m_1,\ldots,m_r,m_1',\ldots,m_r'
\in M$ define the functions $-\varphi_i$ on
$\sigma$ and $\sigma'$ respectively. Identifying
$\Lambda_v$ with $N/(N\cap Span(v_1,\ldots,v_r))$,
parallel transport $T_{\gamma}:\Lambda_v\rightarrow\Lambda_v$ around 
$\gamma$ is given by
$$T_{\gamma}(n)=
n+\sum_{j=1}^r\langle m_j'-m_j,n\rangle (v_j'-v_j).$$
\end{proposition}

Proof. We can identify the tangent space to a point in $\Int(\msig)$
with $Span(m_1,\ldots,m_r)^{\perp}$ and the tangent space to a point
in $\Int(\msig')$ with $Span(m_1',\ldots,m_r')^{\perp}$. 
Then $T_{\gamma}$ is given by the following
chain of identifications in a clockwise order:
$$\begin{matrix}
N/(N\cap Span(v_1,\ldots,v_r))&\mapleft{\psi_1}&
Span(m_1,\ldots,m_r)^{\perp}\\
\mapup{\psi_4}&&\mapdown{\psi_2}\\
Span(m_1',\ldots,m_r')^{\perp}&\mapright{\psi_3}&
N/(N\cap Span(v_1',\ldots,v_r'))\end{matrix}$$
where each arrow is given by projection.
Noting that 
$$\langle m_i,v_j\rangle=\langle m_i',v_j\rangle
=\langle m_i,v_j'\rangle=\langle m_i',v_j'\rangle=-\delta_{ij},$$
we get for $n\in N/(N\cap Span(v_1,\ldots,v_r))$ that 
\begin{eqnarray*}
T_{\gamma}(n)&
=&\psi_4(\psi_3^{-1}(\psi_2(\psi_1^{-1}(n))))\\
&=&\psi_4(\psi_3^{-1}(\psi_2(n+\sum_{j=1}^r\langle m_j,n\rangle v_j)))\\
&=&\psi_4\left(n+\sum_{j=1}^r\langle m_j,n\rangle v_j
+\sum_{k=1}^r\langle m_k', n+\sum_{j=1}^r \langle m_j,n\rangle v_j\rangle
v_k'\right)\\
&=&n+\sum_{k=1}^r (\langle m_k',n\rangle-\langle m_k,n\rangle) v_k'
\mod Span(v_1,\ldots,v_r)\\
&=&n+\sum_{j=1}^r\langle m_j'-m_j\rangle (v_j'-v_j)\mod Span(v_1,\ldots,v_r)
\end{eqnarray*}
as desired.
\qed

\begin{definition}
It follows from Lemma \ref{fan combo} that if $\mtau\in\P_{\nabla}$, then the
fan $\Sigma_{\mtau}$ (see \cite{AffI}, Definition 1.35) 
coincides with the quotient fan
$\Sigma(C(\tau^*))$. Define $\varphi_{\mtau}$ on the fan 
$\Sigma(C(\tau^*))$ by choosing $m\in M$ such that $\langle m,n\rangle
=\varphi(n)$ for all $n\in\tau^*$, so that $\varphi -m$ descends
to a piecewise linear function $\varphi_{\mtau}$ on the fan
$\Sigma(C(\tau^*))=\Sigma_{\mtau}$; the function $\varphi_{\mtau}$ is
well-defined up to a linear function. Such a function also pulls back
to a function on $W_{\mtau}\subseteq B_{\nabla}$ under the projection
$W_{\mtau}\rightarrow N_{\RR}/\RR C(\tau^*)$, 
where $W_{\mtau}$
is the union of interiors of all simplices of $\Bar(\P_{\nabla})$
containing the barycentre of $\mtau$. Then $\varphi_{\nabla}
=\{(W_{\!\mtau},\varphi_{\!\mtau})\}$ defines a multi-valued piecewise
linear function on $B_{\nabla}$. (\cite{AffI}, Definition 1.42.)
Similarly, we obtain $\check\varphi_{\Delta}$
on $B_{\Delta}$ from $\check\varphi$. 
\end{definition}

\begin{theorem} 
$(B_{\Delta},\P_{\Delta},\check\varphi_{\Delta})$ is the discrete
Legendre transform of $(B_{\nabla},\P_{\nabla},\varphi_{\nabla})$.
\end{theorem}

Proof. 
We already have a one-to-one inclusion reversing correspondence between
$\P_{\nabla}$ and $\P_{\Delta}$ by Corollary \ref{correspondence}. To
show the result, we need to show (1) for $\msig\in\P_{\nabla}$ a vertex,
$\sigma$ is the Newton polytope of $\varphi_{\msig}$ on the fan
$\Sigma_{\msig}$, and dually (2) for $\sigma\in\P_{\Delta}$ a
vertex, $\msig$ is the Newton polytope of $\check\varphi_{\sigma}$
on the fan $\Sigma_{\sigma}$.

For (1), $\Sigma_{\msig}=\Sigma(C(\sigma^*))$ by Lemma \ref{fan combo},
and as $\Delta$ is the Newton polytope of the function $\varphi$ on the
fan $\Sigma$, the face $\sigma$ is the Newton polytope of the
function $\varphi_{\msig}$ on $\Sigma_{\sigma}$. For (2) the dual argument
works.
\qed
\bigskip

We now relate these combinatorial constructions to toric degenerations.
Each line bundle $\shL_i$ comes with a canonical section $s_i^0$ given by
$0\in\Delta_i$, whose zero locus is the toric divisor defined by 
$\varphi_i$. Let $s_0:=(s_1^0,\ldots,s_r^0)\in \Gamma(\shL_1\oplus\cdots\oplus
\shL_r)$, and let $s=(s_1,\ldots,s_r)$ be a general section of $\shL_1
\oplus\cdots\oplus\shL_r$. Then
$$ts+s_0=0$$
defines a family $\X$ in $\PP_{\Delta}\times\shS$ with
$\shS=\Spec R$, and $R$ is a discrete valuation ring with
uniformizing parameter $t$. The projection $f:\X
\rightarrow\shS$ gives our family.

Let us assume that the zero locus of $s$ is non-empty and irreducible,
and hence by Proposition \ref{general}, if $s$ is chosen generically, this zero
locus has only canonical singularities.

\begin{proposition}
\label{toricdegen1}
 $f:\X\rightarrow\shS$ is a toric 
degeneration of Calabi-Yau varieties, with dual intersection complex
$B_{\nabla}$, $\P_{\nabla}$.
\end{proposition}

Proof. We first check that $\dim\X_0=n-r$.
Consider the
irreducible components of $\X_0$. The divisor given by $s^0_i=0$
is $\sum\varphi_i(e_j) D_j$, where $D_j$ is the (Weil) divisor
of $\PP_{\Delta}$ corresponding to the ray $e_j$ of $\Sigma$, and 
$\Sigma(1)=\{e_1,\ldots,e_m\}$. Then an irreducible
component of $s_0=0$ is $\bigcap_{j=1}^r D_{i_j}$, 
for some $i_1,\ldots,i_r$ where $\varphi_j(e_{i_j})=1$.
This intersection is either empty, or is a toric
stratum of $\PP_{\Delta}$ defined by the smallest cone $C(\sigma^*)$
of $\Sigma$ containing $e_{i_1},\ldots,e_{i_r}$. Suppose $\dim C(\sigma^*)<r$.
Then $e_{i_1},\ldots,e_{i_r}$ are linearly dependent. However, each
$\varphi_j$ is linear on $C(\sigma^*)$, and since $\varphi_j(e_{i_k})
=\delta_{jk}$, $e_{i_1},\ldots,e_{i_r}$ must be linearly independent.
So $\dim C(\sigma^*)\ge r$, and $\X_0$ has no irreducible component of
dimension $> n-r$, as desired.

Now $\X_{\eta}$ is an irreducible normal variety with canonical singularities
of dimension $n-r$ by Proposition \ref{general}, and $\dim\X=n-r+1$. However, 
$\PP_{\Delta}\times\shS$ is Gorenstein 
and since $\X$ is defined
by $r$ equations, $\X$
can have no components (embedded or otherwise) of dimension $\le n-r$ by
the unmixedness theorem \cite{Matsumura}, Thereom 32.
 Thus $\X$ is irreducible, as $\X_0$ has no component of
dimension greater than $n-r$. Finally, as $\X$ is a complete
intersection in $\PP_{\Delta}\times\shS$, it is also Gorenstein.
Thus $f:\X\rightarrow\shS$ is proper and flat over $S$, and (1) of
Definition \ref{toric degen} holds.

Next we define
$$Z_0=Sing(\X_0)\cap \bigcup_{i=1}^r \{s_i=0\}.$$
Since the $\shL_i$ are semi-ample, and $s$ is general, $s_i$ does
not vanish identically on any toric stratum. 

Let $\bar x\to\X_0$ be a geometric point where none of
the $s_i$'s vanish, with image a point $x\in\X_0$. Then the closure
of the torus orbit of $x$ in $\PP_{\Delta}$ is a toric stratum corresponding
to some cone $C(\sigma^*)$ of $\Sigma$. This cone $C(\sigma^*)$ induces an
affine open subset
$$U_{\sigma}=\Spec k[\dual{C(\sigma^*)}\cap M]$$
of $\PP_{\Delta}$, with $x\in U_{\sigma}$. If $m\in M$, we denote by 
$z^m$ the corresponding monomial.
On $U_{\sigma}$, $\shL_i$ is
trivial, and we can view the sections $s_i$ and $s_i^0$ as regular
functions on $U_{\sigma}$.
Explicitly, if $m_i\in M$ defines the function $\varphi_i$ on the cone
$C(\sigma^*)$, then $m_i\in\dual{C(\sigma^*)}\cap M$ can be taken to
represent the sections $s^0_i$, up to an invertible function on $U_{\sigma}$.
We can then write the equations of $\X$ in $U_{\sigma}\times\shS$ as
$$tf_i-z^{m_i}=0,\quad i=1,\ldots,r$$
with $f_i$ an invertible function times $s_i$. Since 
all the $f_i$'s are non-zero at $x$, a neighbourhood of $x$ is
locally \'etale equivalent to the
subscheme of $U_{\sigma}\times\shS$ defined by
$$t-z^{m_i}=0\quad i=1,\ldots,r.$$
To see this note that if we choose a basis $e_1,\ldots,e_n$ of $M$,
and $g_1,\ldots,g_n$ are functions on $U_{\sigma}\times \shS$ invertible on
a neighbourhood $V$ of $x$, then we obtain in some neighbourhood of $x$ a local
isomorphism $V\rightarrow U_{\sigma}\times\shS$
via $z^{e_i}\mapsto g_iz^{e_i}$. After taking roots of the functions $f_i$,
i.e. after passing to an \' etale cover, we can find $g_i$'s
such that $z^{m_i}\mapsto f_i z^{m_i}$. Under this map the equation 
$tf_i-z^{m_i}=0$ is taken to $t-z^{m_i}=0$.
Note this is locally the same thing as the subscheme of $U_{\sigma}$ defined by 
$$z^{m_i}=z^{m_j},\quad 1\le i<j\le r.$$
This is a toric subvariety of $U_{\sigma}$, defined by the following cone.
Let $V_{\sigma}\subseteq M_{\RR}$ be the subspace spanned by
$\{m_i-m_j| 1\le i<j\le r\}$. Then by Lemma \ref{fan combo},
$$C(\msig)=C(\sigma^*)\cap V_{\sigma}^{\perp}\subseteq V_{\sigma}^{\perp}$$
defines an $n-r+1$ dimensional toric variety $U_{\msig}$ with 
$\dual{C(\msig)}\subseteq \dual{(V_{\sigma}^{\perp})}=M_{\RR}/V_{\sigma}$,
so that the monomials $z^{m_i}$ all agree on $U_{\msig}$. Then the inclusion
$V_{\sigma}^{\perp}\subseteq N_{\RR}$ induces a closed embedding
of $U_{\msig}$ in $U_{\sigma}$, with image defined 
by the equations $z^{m_i}=z^{m_j}$, and the map $U_{\msig}\rightarrow
\Spec k[t]$ is defined by any of the monomials $z^{m_i}$. Thus there is an
\' etale neighbourhood $U_{\bar x}$ of $x$ in $\X$ and a diagram as desired
$$\begin{matrix}
U_{\bar x}&\mapright{}& U_{\msig}\\
\mapdown{f|_{U_{\bar x}}}&&\mapdown{f_x}\\
\Spec R&\mapright{}&\Spec k[\NN]\end{matrix}$$
where $f_x$ is given by any of the $z^{m_i}$. Now since $\langle m_i,n\rangle
=1$ for any $n\in\msig$, $f_x$ vanishes precisely once on each toric
divisor of $U_{\msig}$. This demonstrates that at $\bar x$ the condition of
(4) of Definition \ref{toric degen} holds, and that $\X_0$ is reduced. Note
that if $x\not\in Sing(\X_0)$, $f$ is then smooth at $x$ and the
condition of (4) holds vacuously. In particular, $\X$ is regular in
codimension one and hence normal. We can then take $Z=Z_0\cup
\overline{Sing(\X_{\eta})}$ where $\X_{\eta}$ is the generic fibre
of $f$. From the above it is clear that $Z$ satisfies the conditions
of (4).

The rest of condition (2) of Definition
\ref{toric degen} follows, and condition (3) is easily checked. 

That $B_{\nabla},\P_{\nabla}$ is the dual intersection
complex of $f:\X\rightarrow\shS$ now follows immediately from the construction
described in \S 1: If $\sigma^*\in\R_{\Delta^*}$ is maximal corresponding
to a zero-dimensional stratum of $\X_0$, then by above the corresponding
maximal cell in the dual intersection complex is $\msig$. If $\sigma^*
\in\R_{\Delta^*}$ corresponds to an irreducible component of 
$\X_0$, then this component is defined by the fan $\Sigma(C(\sigma^*))
=\Sigma_{\msig}$, as desired. \qed

\begin{remark}
$B_{\nabla}$ and $B_{\Delta}$ need not be topological spheres: for
example, products of spheres occur. We will not analyze the topology
here, but point out the following elementary observations: in \cite{BB},
the authors introduce the notion of $\Delta_1,\ldots,\Delta_r$
being $k$-independent (\cite{BB}, Def. 3.1) Then it follows from
the proof of \cite{BB}, Theorem 3.3 that if $\Delta_1,\ldots,\Delta_r$
are $k$-independent, $k\ge 3$, then
$$h^1(\O_{\X_0})=\cdots=h^{k-2}(\O_{\X_0})=0.$$
But by \cite{AffI}, Proposition 2.35, $H^i(\X_0,\O_{\X_0})=H^i(B_{\nabla},\CC)$,
so if $\Delta_1,\ldots,\Delta_r$ are $\dim M_{\RR}-r+1$-independent,
then $B_{\nabla}$ is a rational homology sphere. After this paper
was completed, Haase and Zharkov \cite{HZ2} showed that $B_{\nabla}$
is in fact a sphere in this last case.
\end{remark}

\section{The Batyrev-Borisov construction: the general case.}

We now proceed to the general case. We would like to allow a more
general polarization and degeneration.
We are initially given data as in \S 2: $\Delta, \Delta_i, \nabla, \nabla_i, 
\Delta^*,\nabla^*,\varphi=\sum_{i=1}^r\varphi_i$ and $\check\varphi=\sum_{i=1}^r
\check\varphi_i$,
where $\varphi,\varphi_i$ are piecewise linear functions on the
normal fan $\Sigma$ to $\Delta$, 
and $\check\varphi$, $\check\varphi_i$ are piecewise linear functions
on the normal fan $\check\Sigma$ to $\nabla$. 

We will now consider some additional data: subdivisions
$\Sigma'$ and $\check\Sigma'$ of $\Sigma$ and $\check\Sigma$ respectively,
and functions
\begin{eqnarray*}
h:&N_{\RR}\rightarrow\RR\\
\check h:&M_{\RR}\rightarrow\RR
\end{eqnarray*}
which are integral 
piecewise linear and strictly convex on the fans $\Sigma'$ and
$\check\Sigma'$ respectively. These define ample
divisors on partial resolutions of $\PP_{\Delta}$ and $\PP_{\nabla}$
respectively. We will make an additional assumption
that the functions
\begin{eqnarray*}
h':=h-\varphi:&N_{\RR}\rightarrow\RR\\
\check h':=\check h-\check\varphi:&M_{\RR}\rightarrow\RR
\end{eqnarray*}
are convex on $\Sigma'$ and $\check\Sigma'$ respectively (though not necessarily
strictly convex). 

In general, if $g:N_{\RR}\rightarrow\RR$ is a piecewise linear convex
function on $\Sigma'$, write $\Delta^g$ for the Newton polytope of $g$, i.e.
\[
\Delta^g=\{m\in M_{\RR}| \langle m,n\rangle \ge -g(n)\quad\forall n\in
N_{\RR}\},
\]
and similarly if $\check g:M_{\RR}\rightarrow\RR$ is a piecewise linear
convex function on $\check\Sigma'$, write
\[
\nabla^{\check g}
=\{n\in N_{\RR}| \langle m,n\rangle \ge -\check g(m)\quad\forall m\in
M_{\RR}\}.
\]
For any subset $\tau\subseteq\Delta$ for $\Delta\subseteq M_{\RR}$ a
polytope, we write
\[
N_{\Delta}(\tau):=\{n\in N_{\RR}|\hbox{$n|_{\tau}$ is constant
and $\langle m,n\rangle\ge\langle m',n\rangle$ for $m\in\Delta$, $m'\in\tau$}
\},
\]
the normal cone of $\Delta$ along the subset $\tau$.

Let $\Sigma'_{h'}$ be the fan of not necessarily strictly convex cones
on which $h'$ is strictly convex. For example, if $h'$ is
itself linear, $\Sigma'_{h'}$ consists only of the cone $N_{\RR}$.
Similarly define $\check\Sigma'_{\check h'}$.

Cones of $\check\Sigma'$ are in one-to-one correspondence with faces of
$\nabla^{\check h}$, and cones of $\check\Sigma'_{\check h'}$ are in one-to-one
correspondence with faces of $\nabla^{\check h'}$. Write this correspondence
as
\[
\check\delta:\check\Sigma'\rightarrow \{\hbox{faces of $\nabla^{\check h}$}\}
\]
and
\[
\check\delta_{\check h'}:\check\Sigma'_{\check h'}\rightarrow \{\hbox{faces of $\nabla^{\check h'}$}\}.
\]
Generalising the notions of \S 2, if $\sigma\subseteq \partial\Delta^*$
is any set, write
\begin{eqnarray*}
\beta_i^*(\sigma)&=&\{n\in\sigma|\varphi_i(n)=1\}\subseteq\nabla_i\\
\mbeta(\sigma)&=&\beta_1^*(\sigma)+\cdots+\beta_r^*(\sigma)\subseteq
\nabla
\end{eqnarray*}
and for $\sigma\subseteq\partial\nabla^*$, 
\begin{eqnarray*}
\mbeta_i^*(\sigma)&=&\{m\in\sigma|\check\varphi_i(m)=1\}\subseteq\Delta_i\\
\beta(\sigma)&=&\mbeta_1^*(\sigma)+\cdots+\mbeta_r^*(\sigma)\subseteq
\Delta.
\end{eqnarray*}

The faces of a Minkowski sum $\Delta+\Delta'$ are of the form 
$\sigma+\sigma'$ where $\sigma$ and $\sigma'$ are faces of $\Delta$ and
$\Delta'$ respectively, and this decomposition is unique. Then let
\begin{eqnarray*}
\R^{h}_{\nabla^*}&:=&\{(\sigma,\tau)|
\hbox{$\sigma$ is a face of $\nabla^*$ with $\beta(\sigma)\not=\emptyset$,
$\tau$ a face of $\Delta^{h'}$, $\sigma+\tau$ a face of
$\nabla^*+\Delta^{h'}$}\},\\
\R^{\check h}_{\Delta^*}&:=&\{(\sigma,\tau)|
\hbox{$\sigma$ is a face of $\Delta^*$ with $\mbeta(\sigma)\not=\emptyset$,
$\tau$ a face of $\nabla^{\check h'}$, $\sigma+\tau$ a face of
$\Delta^*+\nabla^{\check h'}$}\}.
\end{eqnarray*}
We then define the underlying topological manifolds which will later 
acquire affine structures by
\begin{eqnarray*} 
B^{\check h}_{\nabla}&=&\bigcup_{(\sigma,\tau)\in\R^{\check h}_{\Delta^*}}
\mbeta(\sigma)+\tau,\\
B^h_{\Delta}&=&
\bigcup_{(\sigma,\tau)\in\R^{h}_{\nabla^*}}
\beta(\sigma)+\tau.
\end{eqnarray*}

\begin{lemma}
\label{corrFh}
There is a one-to-one correspondence between proper faces of
$\Delta^*+\nabla^{\check h'}$ and elements of the set
\[
\{
\sigma\cap\tau|\hbox{$\sigma$ is a face of $\Delta$, 
$\tau\in \check\Sigma'_{\check h'},$ $\sigma\cap\tau\not=\emptyset$}\}.
\]
Furthermore, if $\rho$ is an element of this set and
$\rho=\sigma\cap\tau$ for $\sigma$ the minimal face of $\Delta$ containing
$\rho$ and $\tau$ the minimal cone of $\check\Sigma'_{\check h'}$
containing $\rho$, then the corresponding face of $\Delta^*+\nabla^{\check
h'}$ is $\sigma^*+\check\delta_{\check h'}(\tau)$.
\end{lemma}

Proof. Let $\psi:M_{\RR}\rightarrow\RR$ be defined by
\[
\psi(m)=-\inf\{\langle m,n \rangle| n\in\Delta^*\},
\]
so that $\Delta^*=\{n\in N_{\RR}|\hbox{$\langle m,n\rangle \ge -\psi(m)$
for all $m\in M_{\RR}$}\}$. Then by \cite{Oda}, Theorem A.18,
\[
\Delta^*+\nabla^{\check h'}=\{n\in N_{\RR}|\hbox{$\langle m,n\rangle \ge
-\psi(m)-\check h'(m)$ for all $m\in M_{\RR}$}\}.
\]
Furthermore, it is clear that the coarsest decomposition of $M_{\RR}$
into cones on which $\psi+\check h'$ is linear is
$\{C(\sigma)\cap\tau|\hbox{$\sigma$ is a face of $\Delta$, $\tau\in
\check\Sigma'_{\check h'}$}\}$, and by \cite{Oda}, Corollary A.19,
there is a one-to-one correspondence between faces of $\Delta^*
+\nabla^{\check h'}$ with this set of cones. If we only consider
proper faces, then we omit the cone $\{0\}$, and thus the set of
proper faces is in one-to-one correspondence with the given set.
Explicitly, the face of $\Delta^*+\nabla^{\check h'}$ corresponding to
$\sigma\cap\tau$  is
\[
\{n\in N_{\RR}|\hbox{$\langle m,n\rangle \ge -\psi(m)-\check h'(m)$
for all $m\in M_{\RR}$ with equality when $m\in C(\sigma)\cap\tau$}\}.
\]
The last statement then follows easily from this. \qed

\begin{lemma} If $\sigma^*+\check\delta_{\check h'}(\tau)$ is
a face of $\Delta^*+\nabla^{\check h'}$ with $\msig\not=\emptyset$,
then $\msig+\check\delta_{\check h'}(\tau)$ is a face of
$\nabla+\nabla^{\check h'}=\nabla^{\check h}$.
\end{lemma}

Proof. The faces of $\nabla^{\check h}$ are in one-to-one correspondence
with cones $\rho$ of $\check\Sigma'$. If $\rho\in\check\Sigma'$ is given
by $\rho=C(\msig^*)\cap\tau$ for $\msig^*$ a face of $\nabla^*$
and $\tau\in\check\Sigma'_{\check h'}$, and if $C(\msig^*)$
and $\tau$ are minimal cones of $\check\Sigma$ and $\check\Sigma'_{\check h'}$
over $C(\msig^*)\cap\tau$, then the corresponding face of $\nabla^{\check h}$
is $\msig+\check\delta_{\check h'}(\tau)$. Now assume $\sigma^*
+\check\delta_{\check h'}(\tau)$ is a face of $\Delta^*+\nabla^{\check h'}$.
Then $C(\sigma),\tau$ are minimal over $C(\sigma)\cap\tau$. Since
$C(\sigma)\cap\tau\subseteq C(\msig^*)\cap\tau$, $\tau$ must be
minimal over $C(\msig^*)\cap\tau$. If there is a smaller face
$\momega^*\subseteq\msig^*$ of $\nabla^*$ with $C(\msig^*)\cap\tau
=C(\momega^*)\cap\tau$, we wish to show $\sigma\cap\tau=\omega\cap\tau$,
therefore violating the assumption that $\sigma$ was minimal over 
$\sigma\cap\tau$. This follows from 

{\it Claim}. $\sigma\cap\tau=\{m\in C(\msig^*)\cap\tau|\check\varphi_i(m)=1,
\quad 1\le i\le r\}$.

To prove the claim, note the left hand side is clearly contained in the
right hand side. Conversely, if $m\in C(\msig^*)\cap\tau$,
write $m=\sum_i r_i m_i$ with $m_i$ vertices of $\msig^*$
and $r_i\ge 0$. Since each $m_i$ is an integral point of
$\nabla^*$, $\check\varphi_j(m_i)=0$ or $1$, with $1$ occurring for
precisely one value of $j$, so by reindexing we can write
$m=\sum r_{ij}m_{ij}$ with $\check\varphi_k(m_{ij})=\delta_{ik}$. Thus
if $\check\varphi_i(m)=1$, $\sum_j r_{ij}=1$, and $m=\sum_{i=1}^r
m_i'$ with $m_i'=\sum_j r_{ij}m_{ij}$ points of $\msig^*$ with
$\check\varphi_i(m_i')=1$. Thus $m\in\sigma$. \qed

\begin{corollary}
\label{boundaries}
$B_{\nabla}^{\check h}\subseteq\partial\nabla^{\check h}$ and 
$B_{\Delta}^h\subseteq\partial\Delta^h$. \qed
\end{corollary}

Unlike in \S 2, we won't separate the combinatorial description
of the affine structures on these two manifolds from the geometric
arguments, as these structures depend on some additional choices motivated
by the geometry.

\bigskip

In general, $\Delta$ and $\nabla$ define projective toric varieties
$(\PP_{\Delta},\O_{\PP_{\Delta}}(1))$, $(\PP_{\nabla},\O_{\PP_{\nabla}}(1))$
and $\Delta^h$, $\nabla^{\check h}$ define (partial crepant resolutions)
$\pi:\PP_{\Delta^h}\rightarrow\PP_{\Delta}$ and $\check\pi:
\PP_{\nabla^{\check h}}\rightarrow\PP_{\nabla}$, where
$\PP_{\Delta^h}$, $\PP_{\nabla^{\check h}}$ are toric varieties 
defined by the fans $\Sigma'$ and $\check\Sigma'$ respectively.

We could then, for example,
define a family $\X\rightarrow \shS$ with $\X\subseteq\PP_{\Delta}\times\shS$
with equations
$$\sum_{m\in\Delta_1\cap M} t^{\check h(m)} z^m=\cdots=\sum_{m\in\Delta_r
\cap M} t^{\check h(m)}z^m=0$$
where $m\in\Delta_i\cap M$ gives a monomial $z^m$ on the large torus
orbit, extending to a section of the line bundle $\shL_i$
represented by the function $\varphi_i$.

There are two problems with this. First, we need to resolve $\PP_{\Delta}$
so that we can suitably polarize the family using $h$. Second,
the singularities of $\X$ along the fibre with $t=0$ are exceptionally
bad, and these must be resolved.

To resolve these issues, we proceed as follows. 
Let $\tilde\Delta_i\subseteq \tilde M_{\RR}=M_{\RR}\oplus{\RR}$
be defined as
$$
\tilde\Delta_i=\{(m,l)\in M_{\RR}\oplus {\RR}|
m\in\Delta_i, l\ge \check h'(m)\}.
$$
In other words, $\tilde\Delta_i$ is a polytope
extending infinitely upwards but projecting to $\Delta_i$.
$\tilde\Delta_i$ has two types of faces: the {\it vertical} faces
which project to a polytope contained in $\partial\tilde\Delta_i$
of one smaller dimension, and the {\it lower} faces, which project
homeomorphically to subpolytopes of $\Delta_i$. (We remark
that a similar method was used in \cite{Hu} to obtain examples
of semi-stable degenerations of Calabi-Yau manifolds. To apply that
method here, one would use $\check h$ instead of $\check h'$
in the definition of $\tilde\Delta_i$. We do not use this as
this does not yield a toric degeneration.)

Next, let $\tilde\Delta=\tilde\Delta_1+\cdots+\tilde\Delta_r$. 
The polytope $\tilde\Delta$ again extends infinitely upwards,
and projects to $\Delta$. The faces of $\tilde\Delta$ are again
separated into vertical and lower faces. 
Let $\tilde\Sigma$ be the normal fan to
$\tilde \Delta$, living in $\tilde N_{\RR}=N_{\RR}\oplus{\RR}$.
In other words, for any face $\tilde\sigma$ of $\tilde\Delta$,
we have a cone in $\tilde\Sigma$ of the form $N_{\tilde\Delta}(\tilde\sigma)$.

Denote by $X(\tilde\Sigma)$ the $n+1$-dimensional toric variety defined
by the fan $\tilde\Sigma$.

Let $t$ denote
the monomial function on $X(\tilde\Sigma)$ given by
$(0,1)\in \tilde M=M\oplus\ZZ$. A priori this is a rational function.
We have

\begin{proposition}
\label{basicfamily}
The function $t$ is regular, inducing
a map $f:X(\tilde\Sigma)\rightarrow {\AA}^1$.
Each fibre of $f$ except over $t=0$ 
is isomorphic to $\PP_{\Delta}$.
\end{proposition}

Proof. There are two sorts of one-dimensional cones in $\tilde\Sigma$.
One sort are normal to vertical faces of $\tilde\Delta$,
and these are of the form $(n,0)$ for $n\in\Sigma(1)$. The other sort
are normal to maximal lower faces. Since $(0,l)\in\tilde\Delta$ for $l>0$,
if $(m,r)$ is a primitive generator of $N_{\tilde\Delta}(\tilde\sigma)$
for a lower face $\tilde\sigma$,
we must have $rl>0$, so $r>0$. The function $t$
vanishes to order $r$ along the divisor of $X(\tilde\Sigma)$ corresponding
to $(m,r)$; thus in particular $t$ does not have a pole along any
divisor of $X(\tilde\Sigma)$. From this we conclude $t$ is a
regular function.

It is easy to see that if we intersect the fan $\tilde\Sigma$
with $N_{\RR}$ we obtain the fan $\Sigma$. The second statement then
follows from the exact sequence
$$\exact{N_{\RR}}{\tilde N_{\RR}}{\RR}$$
and corresponding maps of fans of $\Sigma$ in $N_{\RR}$ to
$\tilde\Sigma$ in $\tilde N_{\RR}$, and $\tilde\Sigma$ to the
fan in ${\RR}$ defining ${\AA}^1$. The second map of fans induces the
map $f$, and the first map embeds $\PP_{\Delta}$ in $X(\tilde\Sigma)$
as  a general fibre of the map $f$. \qed

We now need a more detailed description of the fan $\tilde\Sigma$, and hence
a more detailed description of $\tilde\Delta$. We do this as follows.

\begin{proposition}
\label{Deltatwiddle}
$$\tilde\Delta=\{(m,l)|m\in \Delta,l\ge \check h'(m)\}.$$
In particular, there is a one-to-one correspondence between maximal lower faces
of $\tilde\Delta$ and maximal cones $\sigma\in\check\Sigma_{\check h'}$,
with
$$\{(m,\check h'(m))|m\in \sigma\cap\Delta\}$$
the maximal lower face of $\tilde\Delta$ corresponding to 
$\sigma\in\check\Sigma_{\check h'}$.
\end{proposition}

Proof.
It is clear that
$$\tilde\Delta\subseteq \{(m,l)|m\in \Delta,l\ge\check h'(m)\};$$
indeed, if $(m,l)\in\tilde\Delta$, then $(m,l)=\sum_{i=1}^r (m_i,l_i)$
for some $(m_i,l_i)\in\tilde\Delta_i$, so $l_i\ge \check h'(m_i)$.
But by convexity of $\check h'$,
$$\check h'(m)\le\sum \check h'(m_i)\le \sum_{i=1}^r l_i=l.$$

We only need to show the converse inclusion now. To do this, let $\sigma$
be any maximal cone of $\check\Sigma'$, set 
\begin{eqnarray*}
\tau_i&:=& \mbeta_i^*(\sigma\cap\partial\nabla^*)\\
\sigma_i&:=& Conv\{0,\tau_i\}
\end{eqnarray*}
and set $\sigma':=\sigma_1+\cdots+\sigma_r$. 

\begin{claim}
\label{mink1}
$\sigma'=\sigma\cap\Delta$.
\end{claim}

Proof. Since $\Delta=\Delta_1+\cdots+\Delta_r$, it is clear that $\sigma'
\subseteq \sigma\cap\Delta$.
If $m\in\sigma_i$ then $lm\in \sigma_i$
for $0\le l\le 1$. Thus $Conv\{\sigma_1,\ldots,\sigma_r\}\subseteq\sigma_1
+\cdots+\sigma_r$. On the other hand, $\sigma\cap\partial\nabla^*$ is
the convex hull of some integral points; let $v$ be such an integral
point. Then $\check\varphi_i(v)=1$ for some $i$,
so $v\in\sigma_i$ for some $i$. Thus $Conv\{\sigma_1,\ldots,\sigma_r\}
=\sigma\cap\nabla^*$. So we have
\[
\sigma\cap\nabla^*\subseteq\sigma'\subseteq\sigma\cap\Delta.
\]
Now let $\rho$ be a maximal face of $\sigma'$ not contained in
a face of $\sigma$. If we can prove $\rho\subseteq\partial\Delta$, then
we will obtain the equality $\sigma'=\sigma\cap\Delta$.
It is clear 
that $0\not\in\rho$ since $\sigma\cap\nabla^*\subseteq\sigma'$. 
Since $\rho=\rho_1+\cdots+\rho_r$ where $\rho_i$ is a face of $\sigma_i$,
we must have $0\not\in\rho_i$ for at least one $i$.

Since $\sigma$ is maximal, there exists a unique vertex $n$ of 
$\nabla$ such that $\langle n,m\rangle =-1$ for
all $m\in \sigma\cap\partial\nabla^*$, and since $\nabla=\nabla_1+\cdots
+\nabla_r$, $n$ decomposes uniquely as $n=n_1+\cdots+n_r$ with
$n_i$ a vertex of $\nabla_i\subseteq\Delta^*$ and $n_i$ defines
the function $-\check\varphi_i$ on $\sigma$.
Now $\tau_i$ is the unique
maximal face of $\sigma_i$ not containing $0$ (note for some $i$,
$\tau_i$ may be empty, but this is only the case when $\sigma_i=\{0\}$).
Thus $\rho_i\subseteq\tau_i$ so 
$$\rho\subseteq \sigma_1+\cdots+\tau_i+\cdots+\sigma_r.$$
However $n_i\in\Delta^*$, and $\langle n_i,\sigma_j\rangle=0$, so 
$\langle n_i,m\rangle=-1$ for all $m\in\rho$. Thus $\rho$ is contained in
$\partial\Delta$. \qed 

Now for any $m\in\Delta$, there is some maximal $\sigma\in\check\Sigma'$ such that
$m\in\sigma$, and thus
$m\in\sigma'$ by the claim. Then 
$m=\sum m_i$ with $m_i\in\sigma_i$ and $(m,\check h'(m))
=\sum (m_i,\check h'(m_i))\in\tilde\Delta$ as $m_1,\ldots,m_r\in\sigma$
on which $\check h'$ is linear. \qed

\begin{proposition} 
\label{sigmatwiddle}
\begin{eqnarray*}
\tilde\Sigma&=&\{C(\sigma^*)\times\{0\}|\hbox{$\sigma^*$ a face of $\Delta^*$}\}
\\
&\cup& \left\{C(\sigma^*)\times\{0\}+C(\tau\times\{1\})\bigg|
{\hbox{$\sigma^*$ a face of $\Delta^*$, $\tau$ a face of $\nabla^{\check h'}$,}
\atop\hbox{
and
$\sigma^*+\tau$ a face of $\Delta^*+\nabla^{\check h'}$}}\right\}\\
&\cup&\{C(\tau\times\{1\})|
\hbox{$\tau$ a face of $\nabla^{\check h'}$}\}.
\end{eqnarray*}
Here $C$ denotes the cone over a set with vertex the origin, as usual. 
\end{proposition}

Proof. The normal cone to a vertical face of $\tilde\Delta$ mapping to a face
$\sigma$ of $\Delta$ is just $C(\sigma^*)\times\{0\}$. A lower
face $\rho$ of $\tilde\Delta$ is projected either to a subpolytope
of $\Delta$ containing $0$, or is mapped to a polytope contained in 
$\partial\Delta$. In the former case, $\rho$ maps to $\tau\cap\Delta$
for some $\tau\in\check\Sigma'_{\check h'}$. The maximal faces
of $\tilde\Delta$ containing $\rho$ are then in one-to-one correspondence
with maximal cones $\tau'$ of $\check\Sigma'_{\check h'}$ containing
$\tau$, and then the normal cone to $\tilde\Delta$ along $\rho$ is generated
by the set of normal vectors to the maximal faces of $\tilde\Delta$
containing $\rho$, and this is
\[
\{(-n'_{\tau'},1)|\hbox{$\tau'\in\Sigma'_{\check h'}$ maximal cones containing
$\tau$}\}
\]
where $n'_{\tau'}$ defines $\check h'$ on $\tau'$.
But $\check\delta_{\check h'}(\tau)$ is the convex hull of 
\[
\{-n'_{\tau'}|\hbox{$\tau'\in\Sigma'_{\check h'}$ a maximal cone containing
$\tau$}\},
\]
so the normal cone to $\tilde\Delta$ along $\rho$ is 
$C(\check\delta_{\check h'}(\tau)\times\{1\})$. 
In the second case, $\rho$ is the intersection of some minimal vertical
face containing it, mapping, say, to a face $\sigma$ of $\Delta$,
and a minimal horizontal face containing both $\rho$ and $0$,
mapping to $\tau\cap\Delta$ for some $\tau\in\check\Sigma'_{\check h'}$.
In this case, the normal cone is
\[
C(\sigma^*)\times\{0\}+C(\check\delta_{\check h'}(\tau)\times\{1\}).
\]
But the set of such $\rho$ is clearly in one-to-one correspondence with the
set defined in Lemma \ref{corrFh}, so the result follows.
\qed
\bigskip

The polytope $\tilde\Delta$ determines a line bundle $\tilde\shL$ on
$X(\tilde\Sigma)$. This line bundle is induced by the piecewise
linear function $\tilde\varphi$ on $\tilde\Sigma$ defined by
$$\tilde\varphi(\tilde n)=-\inf\{\langle \tilde n,\tilde m\rangle|
\tilde m\in\tilde\Delta\}.$$
We clearly have $\tilde\varphi(n,0)=1$ for $n\in\partial\Delta^*$,
and since $(0,0)$ is in every maximal lower face of $\tilde\Delta$,
$\tilde\varphi(n,1)=0$ for every 
$n\in\nabla^{\check h'}$. Also, the Minkowski
decomposition $\tilde\Delta=\tilde\Delta_1+\cdots+\tilde\Delta_r$ induces
a decomposition $\tilde\shL=\tilde\shL_1\otimes\cdots\otimes\tilde\shL_r$.
This induces a decomposition
$\tilde\varphi=\tilde\varphi_1+\cdots+\tilde\varphi_r$, with 
$\tilde\varphi_i(n,0)=\varphi_i(n)$ and $\tilde\varphi_i(n,1)=0$ for any 
$n\in\nabla^{\check h'}$.

\begin{definition}
\label{goodsub}
 A subdivision $\tilde\Sigma'$ of the fan $\tilde\Sigma$
is {\it good} for $\Sigma'$ if the fan 
$\{\sigma\cap(N_{\RR}\times\{0\})|\sigma\in\tilde\Sigma'\}$ coincides
with $\Sigma'$ and furthermore, every one-dimensional cone
of $\tilde\Sigma'$ not contained in $N_{\RR}\times\{0\}$ is 
generated by a primitive vector $(n,1)$ with 
$n\in \nabla^{\check h'}$. 
\end{definition}

\begin{observation}
\label{conetypes}
If $\tilde\Sigma'$ is a good subdivision for $\Sigma'$, it contains
three sorts of cones:
\begin{enumerate}
\item Cones of the form $C(\sigma)\times\{0\}$ for $C(\sigma)\in
\Sigma'$ with $\sigma\subseteq\partial\Delta^*$;
\item
Cones of the form $C(\sigma)\times\{0\}+C(\tau\times\{1\})$
where $C(\sigma)\in\Sigma'$ for $\sigma\subseteq\partial\Delta^*$,
$\tau\subseteq\partial\nabla^{\check h'}$, $\sigma+\tau$ contained
in a face of $\Delta^*+\nabla^{\check h'}$.
\item Cones contained in $C(\nabla^{\check h'}\times\{1\})$.
\end{enumerate}
We will call a cone of $\tilde \Sigma'$ {\it relevant} if it is of the
second type with $\mbeta(\sigma)\not=\emptyset$.
\end{observation}
\bigskip
Given a subdivision, we obtain a (partial) resolution $\tilde\pi:
X(\tilde\Sigma')\rightarrow X(\tilde\Sigma)$. There is a function
$f':X(\tilde\Sigma')\rightarrow\AA^1$ ($f'=f\circ \tilde\pi$), again coming
from the monomial $t$, whose general fibre, as in Proposition \ref{basicfamily},
is isomorphic to $\PP_{\Delta^h}$.

Let $\tilde\shL_i'=\tilde\pi^*\tilde\shL_i$. The functions 
$\tilde\varphi_i$ are the piecewise linear functions on $\tilde\Sigma'$
determining $\tilde\shL_i'$ also, and again the Newton polytopes of
$\tilde\shL_i'$ are $\tilde\Delta_i$.

Let $s_i\in\Gamma(X(\tilde\Sigma'),\tilde\shL_i')$ be a general section,
and let $s_i^0\in\Gamma(X(\tilde\Sigma'),\tilde\shL_i')$ be the section
defined by $(0,0)\in\tilde\Delta_i$. 
Thus in particular
$ts_i$ is also a section of $\tilde\shL_i$. Let $\X'\subseteq X(\tilde\Sigma')$
be defined by the equations
$$ts_1+s_1^0=\cdots=ts_r+s_r^0=0$$
in $X(\tilde\Sigma')$. We obtain by restriction $f':\X'\rightarrow\AA^1$.
We can replace $\X'$ with a family over $\shS
=\Spec R$, where $R$ is a discrete valuation
ring with uniformizing parameter $t$, via basechange.

\begin{theorem}
\label{toricdegen2}
$f':\X'\rightarrow\shS$ is a toric 
degeneration of Calabi-Yau manifolds.
\end{theorem}

Proof. We follow the same outline as the proof of Proposition \ref{toricdegen1}.
First let us check that $\X'$ is irreducible, $f':\X'\rightarrow\shS$
flat of relative dimension $n-r$. If $\eta$ is the generic point
of $\shS$, then $\X'_{\eta}\subseteq \PP_{\Delta^h}\times\eta$ is cut out
by generic sections of $\pi^*\shL_1,\ldots,\pi^*\shL_r$, and hence is a codimension
$n-r$ complete intersection with canonical singularities. 
We need to check that $\dim\X'_0=n-r$.
Now $\X'_0$ is defined by the equations
$$t=s_1^0=\cdots=s^0_r=0$$
in $X(\tilde\Sigma')$. 
Each of these equations defines a toric divisor, so each irreducible
component of $\X_0'$ is a toric stratum.
Let $\rho$ be a cone in $\tilde\Sigma'$ corresponding to a toric
stratum of $\X_0'$. Then none of $\tilde\varphi_i$ or $(0,1)$ can
be identically zero on the cone $\rho$. Thus following Observation
\ref{conetypes}, $\rho$ cannot be of the first or third type, since
$(0,1)$ vanishes on cones of the first type and $\tilde\varphi_i$
vanishes on cones of the third type. Thus we can write
\[
\rho=C(\sigma)\times\{0\}+C(\tau\times\{1\})
\]
with $\sigma\subseteq\partial\Delta^*$. In addition, we must have
$\mbeta(\sigma)\not=\emptyset$, since for each $i$, $\tilde\varphi_i$
must take the value $1$ on some vertex of $\sigma$. Thus $\rho$ is
a relevant cone.

Now suppose $\rho$ corresponds to an irreducible component of $\X_0'$.
Then $\sigma$ has vertices $n_1,\ldots,n_r$ with $\varphi_i(n_j)
=\delta_{ij}$. Thus as in the proof of Proposition \ref{toricdegen1},
$n_1,\ldots,n_r$ are linearly independent, so $\dim\rho\ge r+1$.
Thus the dimension of the corresponding
irreducible component is $\le n-r$. Again by the unmixedness theorem, it
must be dimension $n-r$. (Any toric variety is Cohen-Macaulay.)
Thus we can conclude
as in Proposition \ref{toricdegen1} that $\X'$ is irreducible.
Thus $f':\X'\rightarrow \shS$ is proper and flat, and $\X'$ is
irreducible, and (1) of Definition \ref{toric degen} holds.

Next we define the set $Z_0$ as
$$Z_0=Sing(\X_0')\cap
\bigcup_{i=1}^r\{s_i=0\}.$$
Since the $s_i$ are general, they do not vanish identically on any toric 
stratum. 

We will now consider (4) of Definition \ref{toric degen}. Let
$\bar x\rightarrow\X_0'$ be a geometric point in $\X_0'$ where
none of the $s_i$'s vanish, with image $x\in\X_0'$.
Then the closure of the torus orbit of $x$ in $X(\tilde\Sigma')$
is a toric stratum corresponding to some cone $\rho$ of $\tilde\Sigma'$. 
This cone $\rho$ induces an affine open subset
$$U_{\rho}=\Spec k[\dual{\rho}\cap (M\oplus\ZZ)]$$
of $X(\tilde\Sigma')$. On $U_{\rho}$, we can view $s_i$ and $s_i^0$
as regular functions.
Explicitly, if $(m_i,l_i)\in
M_{\RR}\oplus\RR$ are chosen to define the function $\tilde\varphi_i$
on the cone $\rho$ (with $m_i$ defining $\varphi_i$ on $\rho
\cap N_{\RR}$) then thinking of $z^{(m_i,l_i)}$ as a monomial function
on $U_{\rho}$, we can write the equation of $\X'$ in $U_{\rho}$ as
$$tf_i-z^{(m_i,l_i)}=0,\quad i=1,\ldots,r$$
with $f_i$ proportional to $s_i$. Since all the $f_i$'s are non-zero
at $x$, a neighbourhood of $x$ in $\X'$ 
is locally \'etale equivalent (as in the proof of Proposition 
\ref{toricdegen1})
to the subscheme of $U_{\rho}$ defined by
$$t-z^{(m_i,l_i)}=0,\quad i=1,\ldots,r.$$
This is a toric subvariety $U_{\rho'}$ of $U_{\rho}$, defined by the 
following cone. Let $V_{\rho}\subseteq M_{\RR}\oplus\RR$ be
the subspace spanned by $\{(m_i,l_i-1)|1\le i\le r\}$. Then the cone
$\rho':=\rho\cap V_{\rho}^{\perp}$ defines an $n-r+1$ dimensional
toric variety $U_{\rho'}$. In addition, the map $f'$ is \'etale
equivalent to the induced map $U_{\sigma'}\rightarrow\Spec k[t]$
defined by $t$.

\begin{claim}
\label{basicclaim}
$\rho'=C((\mbeta(\sigma)+\tau)\times\{1\})$.
\end{claim}

Proof. Clearly
\[
\rho'=\{(n,l)\in\rho|\tilde\varphi_i(n,l)=l\quad\forall i\}.
\]
Write $(n,l)\in\rho$ as $(n,l)=(n',0)+l(n'',1)$ with $n'\in C(\sigma)$
and $n''\in\tau$. Then $(n,l)\in\rho'$ if and only if $\varphi_i(n')=l$ for
all $i$.
Now suppose $n'$ is a vertex of $\mbeta(\sigma)$. Then $n'=n_1+\cdots+n_r$,
with $n_i\in\beta_i^*(\sigma)$, so $\varphi_i(n_j)=\delta_{ij}$. Thus
$\varphi_i(n')=1$, so $(n',0)+(n'',1)\in\rho'$ for any $n''\in\tau$.
Thus $C((\mbeta(\sigma)+\tau)\times\{1\})\subseteq\rho'$. Conversely,
suppose $(n,l)\in\rho'$. If $l=0$, then $\varphi_i(n)=0$ for all $i$,
so $\varphi(n)=0$. But as $\varphi$ takes the constant value $1$
on $\partial\Delta^*$, this implies $n=0$. Thus we can assume $l\not=0$,
and by dividing by $l$, assume $(n,1)\in\rho'$. Of course,
we have $(n,1)=(n',0)+(n'',1)$ as above, and we just need to show
$n'\in\mbeta(\sigma)$. Write $n'=\sum r_jn_j$ for $n_j$ vertices of
$\sigma$, $r_j>0$.  Then as each $n_j$ is integral, we have $\varphi_i(n_j)
=0$ or $1$, and since $\varphi=\sum\varphi_i$, $\varphi_i(n_j)=1$ for some
$i$. Thus, reindexing, we can write $n'=\sum_{(i,j)\in I} r_{ij}n_{ij}$
with $\varphi_k(n_{ij})=\delta_{ik}$. Thus $\varphi_i(n')=1$ implies 
$\sum_j r_{ij}=1$, so $\sum_j r_{ij} n_{ij}\in \beta_i^*(\sigma)$.
Thus $n'\in \mbeta(\sigma)$ as desired. \qed

To conclude the proof of Condition (4) of Definition \ref{toric degen},
we just have to observe that from the claim, the function $t$ vanishes
on every toric divisor of $U_{\rho'}$,
in fact to order $1$, so in addition $\X_0$ is reduced. Then 
taking $Z=Z_0\cup \overline{Sing(\X'_{\eta})}$, $Z$ satisfies 
Condition (4).

Finally, to check (2) and (3), the key point is to make sure that for any
irreducible component $X$ of $\X_0'$, and toric divisor $D\subseteq X$,
there exists another irreducible component $Y$ of $\X_0'$ with $X
\cap Y=D$. However, this is easily verified.
\qed

We can now identify the relevant polyhedral decomposition.

\begin{proposition}
\[
B^{\check h}_{\nabla}=\bigcup \mbeta(\sigma)+\tau,
\]
where the union is over all relevant cones $C(\sigma)\times\{0\}+
C(\tau\times\{1\})$ of a subdivision $\tilde\Sigma'$ of $\tilde\Sigma$ 
good for $\Sigma'$.
\end{proposition}

Proof. The inclusion of the right-hand set in the left-hand set is clear. Conversely,
let $n\in B_{\nabla}^{\check h}$. Then $n\in\mbeta(\sigma)+\tau$ for a cone
$C(\sigma)\times\{0\}+C(\tau\times\{1\})$ of $\tilde\Sigma$, and thus $(n,1)$ is in a cone
$C(\sigma')\times\{0\}+C(\tau'\times\{1\})$ of $\tilde\Sigma'$. 
We only need to show $n\in\mbeta(\sigma')
+\tau'$. Certainly $(n,1)=(n',0)+(n'',1)$ for $n'\in C(\sigma')$, $n''\in\tau'$.
Also, it follows from Claim \ref{basicclaim}
applied to $\tilde\Sigma$ rather than $\tilde\Sigma'$, that 
$\tilde\varphi_i(n,1)=1$ for all $i$, so $\varphi_i(n')=1$ for all $i$.
But we can then conclude, just as in the proof of the claim, that $n'\in
\mbeta(\sigma')$. \qed

\begin{definition} We set
\[
\P_{\nabla}^{\tilde\Sigma'}:=\{\mbeta(\sigma)+\tau
|\hbox{$C(\sigma)\times\{0\}+C(\tau\times\{1\})$ is a relevant
cone of $\tilde\Sigma'$}\}.
\]
This is clearly a polyhedral decomposition of $B^{\check h}_{\nabla}$.
\end{definition}

We can then define as usual $\Gamma_{\nabla}^{\tilde\Sigma'}\subseteq
B_{\nabla}^{\check h}$ the union of all simplices of $\Bar(\P_{\nabla}^{\tilde
\Sigma'})$ not containing a vertex of $\P^{\tilde\Sigma'}_{\nabla}$
or intersecting the interior of a maximal face of $\P_{\nabla}^{\tilde\Sigma'}$.
We define an integral affine structure on $B_{\nabla}^{\check h}\setminus
\Gamma_{\nabla}^{\tilde\Sigma'}$ as follows.

For every maximal face $\sigma\in\P^{\tilde\Sigma'}_{\nabla}$ of $B$, 
let ${\AA}_{\sigma}$ be the 
affine subspace of $N_{\RR}$ spanned by $\sigma$. Then we define an
affine chart on $\Int(\sigma)$ by $\psi_{\sigma}:\Int(\sigma)\hookrightarrow
{\AA}_{\sigma}$.

Secondly, consider a vertex $v$ of $\P^{\tilde\Sigma'}_{\nabla}$. 
Then $v$ arises from a relevant cone $C(\sigma)\times\{0\}+C(\tau\times\{1\})$
with $\{v\}=\mbeta(\sigma)+\tau$. Then clearly $\dim\tau=\dim\mbeta(\sigma)=0$,
so $\tau=\{v'\}$ and $\sigma=Conv\{n_1,\ldots,n_r\}$ where
$\varphi_i(n_j)=\delta_{ij}$. Clearly $n_1,\ldots,n_r$ are linearly independent.
Letting $W_v$ be the union of interiors of all simplices of $\Bar(\P_{\nabla}^{\tilde
\Sigma'})$ containing $v$, we define a chart
$$\psi_v:W_v\rightarrow N_{\RR}/Span(n_1,\ldots,n_r)$$
via projection.

\begin{proposition} The charts $\psi_v$ and $\psi_{\sigma}$
define an integral affine structure on $B_{\nabla}^{\check h}\setminus
\Gamma_{\nabla}^{\tilde\Sigma'}$,
making $B_{\nabla}^{\check h}$ an integral affine manifold with singularities
of dimension $\dim M_{\RR}-r$. In addition, this structure makes
$(B_{\nabla}^{\check h},\P_{\nabla}^{\tilde\Sigma'})$ the dual intersection
complex of the degeneration $f':\X'\rightarrow\shS$.
\end{proposition}

Proof.  Let $v$ be a vertex of $\P_{\nabla}^{\tilde\Sigma'}$ 
corresponding to a cone 
\[
\omega_v=C(Conv\{n_1,\ldots,n_r\})\times\{0\}
+C(\{w\}\times\{1\})
\]
of $\tilde\Sigma'$, with $v=n_1+\cdots+n_r+w$, and let $\rho$ be
a maximal facet of $\P_{\nabla}^{\tilde\Sigma'}$ containing $v$,
corresponding to a maximal cone
\[
\omega_{\rho}=C(\sigma)\times\{0\}+C(\tau\times\{1\})
\]
of $\tilde\Sigma'$ with $\rho=\mbeta(\sigma)+\tau$. Then
by Claim \ref{basicclaim}, 
\[
{\AA}_{\sigma}=\{n\in N_{\RR}| (n,1)\in V_{\rho}^{\perp}\}.
\]
More explicitly, if $(m_i,l_i)$ defines $\tilde\varphi_i$ on
$\omega_\rho$, then $-(m_i,l_i)$ is a vertex of $\tilde\Delta_i$, and
in particular is in a lower face of $\tilde\Delta_i$, so $l_i=
-\check h'(m_i)$. Then 
\[
\AA_{\sigma}=
\{n\in N_{\RR}|\langle m_i,n\rangle=\check h(m_i),
i=1,\ldots,r\}.
\]
Now $\langle m_i,n_j\rangle=\delta_{ij}$,
so if $\sum a_jn_j\in {\AA}_{\sigma}$, we have
$a_i=\langle m_i,\sum a_jn_j\rangle=\check h(m_i)$, so
$${\AA}_{\sigma}\cap Span(n_1,\ldots,n_r)=\left\{\sum \check h(m_i)m_i\right
\}.$$
In particular, the projection
${\AA}_{\sigma}\rightarrow N_{\RR}/Span(n_1,\ldots,n_r)$ is an affine
isomorphism, with the linear part preserving the integral structure.

We then have to show that $\psi_v$ is one-to-one. First let's look at
the tangent wedge to $\psi_v(\rho)$ at $\psi_v(v)\equiv w
\mod Span(n_1,\ldots,n_r)$.
As in \S 2, it is easy to see
that the tangent wedge of $\psi_v(\rho)$ at $\psi_v(v)$ is the cone
$(\omega_{\rho}+{\RR}\omega_v)/{\RR}\omega_v$ of
the fan $\tilde\Sigma'(\omega_v)$. Since the maximal facets of
$\P_{\nabla}^{\tilde\Sigma'}$ 
containing $v$ are in one-to-one correspondence with
maximal cones of $\tilde\Sigma'$ containing $\omega_v$,
we see that the tangent wedges to all the $\psi_v(\rho)$'s fit together
to form the fan $\tilde\Sigma'(\omega_v)$. In particular $\psi_v$
is a ($C^0$) immersion.

It is now easy to see $B^{\check h}_{\nabla},\P_{\nabla}^{\tilde\Sigma'}$ 
is the dual intersection complex of the toric
degeneration $f':\X'\rightarrow S$ as in Proposition \ref{toricdegen1}. 
\qed

We compute monodromy.
\begin{proposition} 
\label{monodromy2}
Let $v$ and $v'$ be two vertices of $\P_{\nabla}^{\tilde
\Sigma'}$
and let $\rho$ and $\rho'$ be two
maximal facets of $\P_{\nabla}^{\tilde\Sigma'}$
containing $v$ and $v'$ coming from cones $\omega_{\rho}$ and $\omega_{\rho'}$
of $\tilde\Sigma'$ respectively. Let $\gamma$ be
a simple loop based at $v$, passing successively into $\Int(\rho)$,
through $v'$, into $\Int(\rho')$ and back to
$v$. Write $v=n_1+\cdots+n_r+w$ and $v'=n_1'+\cdots+n_r'+w'$,
and let $(m_1,l_1),\ldots,(m_r,l_r),(m_1',l_1'),
\ldots,(m_r',l_r')
\in M$ define the functions $-\tilde\varphi_i$ on
$\omega_{\rho}$ and $\omega_{\rho'}$ respectively. Identifying
$\Lambda_v$ with $N/(N\cap Span(n_1,\ldots,n_r))$, then
parallel transport $T_{\gamma}:\Lambda_v\rightarrow\Lambda_v$ around 
$\gamma$ is given by
$$T_{\gamma}(n)=
n+\sum_{j=1}^r\langle m_j'-m_j,n\rangle (n_j'-n_j).$$
\end{proposition}

Proof. This proof is identical to that of Proposition \ref{monodromy1}
\qed

We now determine when $B$ is simple. 
(See \cite{AffI}, Definition 1.60). We do this both so
we know when Theorem \ref{torusbundle} can be applied and to
justify the relevance of the simplicity condition.
Recall that an MPCP
(maximal projective crepant partial) resolution of $\PP_{\Delta}$
is obtained by choosing $h$ so that the induced
decomposition of $\partial\Delta^*$ consists only of elementary simplices,
i.e. simplices whose only integral points are the vertices.

\begin{theorem}
\label{simplicity}
 If $h$ and $\check h$ induce MPCP resolutions of
$\PP_{\Delta}$ and $\PP_{\nabla}$ respectively, then 
$(B_{\nabla}^{\check h},\P_{\nabla}^{\tilde\Sigma'})$ is simple.
\end{theorem}

Proof. We follow the notation of \cite{AffI}, \S 1.5. Let $\tau\in
\P_{\nabla}^{\tilde\Sigma'}$, $1\le\dim\tau\le n-r-1$,
\begin{eqnarray*}
\P_1(\tau)&=&\{e:\omega\rightarrow\tau|\dim\omega=1\}\\
\P_{n-r-1}(\tau)&=&\{f:\tau\rightarrow\rho|\dim\rho=n-r-1\}.
\end{eqnarray*}
(Note that in our case, using morphisms $e:\omega\rightarrow\tau$,
etc. is superfluous; whenever $\omega\subseteq\tau$ there is a unique morphism
$\omega\rightarrow\tau$, as there are no self-intersections of cells). 
If $e\in\P_1(\tau)$, $e:\omega\rightarrow\tau$, then we can
calculate $\check\Delta_e(\tau)$ using Proposition \ref{monodromy2}. Choose
some $\sigma,\sigma'\in\P_{\nabla}^{\tilde\Sigma'}$ maximal containing
$\tau$, and let $v,v'$ be the two endpoints of $\omega$. Write
$\omega=\mbeta(\omega_1)+\omega_2$ with
$C(\omega_1\times\{0\})+C(\omega_2\times\{1\})\in\tilde\Sigma'$,
$C(\omega_1)\in\Sigma'$. By the MPCP assumption, $\omega_1$ is
simplicial, with generators $n_1,\ldots,n_k$, $k\ge r$. If 
$k>r+1$, then $\dim\mbeta(\omega_1)\ge 2$, contradicting $\dim\omega=1$,
so we have two cases. Either $\omega_1=Conv\{n_1,\ldots,n_r\}$
with $\varphi_i(n_j)=\delta_{ij}$, or $\omega_1=Conv\{n_1,\ldots,n_r,
n_p'\}$ with $\varphi_i(n_j)=\delta_{ij}$ and $\varphi_p(n'_p)=1$
for some $p$. In the former case, it follows from Proposition \ref{monodromy2}
that there is no monodromy around the loop $\gamma$ determined by the data 
$v,v',\sigma,\sigma'$. In the second case,
\[
v=n_1+\cdots+n_r+w,\quad v'=n_1+\cdots+n_p'+\cdots+n_r+w'
\]
for $w,w'\in\omega_2$, and 
\[
T_{\gamma}(n)=n+\langle m_p'-m_p,n\rangle (n_p'-n_p).
\]
Noting again by the MPCP assumption $n_p'-n_p$ is primitive, and it follows
that we can take, up to translation,
\[
\check
\Delta_e(\tau)=Conv\left\{m_p\bigg|{\hbox{$(m_p,l_p)$ defines $-\tilde\varphi_p$
on the cone of $\tilde\Sigma'$}\atop\hbox{corresponding to $\sigma$, $\tau
\subseteq\sigma$, $\sigma$ maximal}}\right\}.
\]
This only depends on the index $p$, which depends only on $e$, so
we write this index as $p_e$. Then we set for $1\le i \le r$,
\[
\Omega_i=\{e:\omega\rightarrow\tau|\dim \omega_1=r+1, p_e=i\}
\]
and for $e\in\Omega_i$,
\[
\check\Delta_e(\tau)=
\delta_i:=Conv\left\{m_i\bigg|{\hbox{$(m_i,l_i)$ defines $-\tilde\varphi_i$
on
the cone of $\tilde\Sigma'$}\atop\hbox{ corresponding to $\sigma$, $\tau
\subseteq\sigma$, $\sigma$ maximal}}\right\}.
\]
Now $\delta_i$ is the image under projection of a face of $\tilde\Delta_i$,
determining the function
$-\tilde\varphi_i$ on the cone of $\tilde\Sigma'$ corresponding
to $\tau$. It is then easy to see that if this cone is
$C(\tau_1\times\{0\})+C(\tau_2\times\{1\})$, and $\delta_i'$ is the
set of elements of $\Delta_i$ giving the value of $-\varphi_i$
on $\tau_1$ and $\delta''=N_{\nabla^{\check h'}}(\tau_2)\in
\check\Sigma'_{\check h'}$, then $\delta_i=\delta_i'\cap\delta''$.
Then by Lemma \ref{corrbetweenfaces}, $\delta'=Conv(\delta_1',\ldots,
\delta_r')$ is a face of $\nabla^*$ and also 
$Conv(\delta_1,\ldots,\delta_r)=\delta'\cap\delta''$, so the cone over
$Conv(\delta_1,\ldots,\delta_r)$ is a cone in $\check\Sigma'$,
hence $Conv(\delta_1,\ldots,\delta_r)$ is an elementary simplex.
From this one concludes that in $M_{\RR}\times\RR^r$,
$Conv(\delta_1\times\{e_1\},\ldots,\delta_r\times\{e_r\})$ is an elementary
simplex with $e_1,\ldots,e_r$ the standard basis of $\RR^r$. 
This is the first condition needed for simplicity. The second
statement is essentially dual. Taking $f\in\P_{n-r-1}(\tau)$, a similar
argument shows that 
\[
\Delta_f(\tau)=\beta_i^*(\tau_1)
\]
for some $i$ and the cone
over $\tau_1=Conv(\beta_1^*(\tau_1),\ldots,\beta_r^*(\tau_1))$ is a cone
in $\Sigma'$, so by the MPCP assumption we conclude
$Conv(\beta_1^*(\tau_1)\times\{e_1\},\ldots,\beta_r^*(\tau_r)\times
\{e_r\})$ is an elementary simplex. \qed
\bigskip

We now wish to discuss duality. To do so, we need to polarize
the degeneration $f:\X'\rightarrow\shS$, and we do this by
suitably polarizing the ambient space $X(\tilde\Sigma')$. We will
use the following:

\begin{definition}
A subdivision $\tilde\Sigma'$ of $\tilde\Sigma$ with an integral piecewise linear
function $\tilde h$ on $\tilde\Sigma'$ is good for $\Sigma',h$ if
$\tilde\Sigma'$ is a good subdivision of $\tilde\Sigma$ (Definition 
\ref{goodsub}), $\tilde h|_{N_{\RR}\times\{0\}}=h$, $\tilde h$ is
strictly convex,
and $\tilde h':=\tilde h-\tilde\varphi$ is convex on $\tilde\Sigma'$.
\end{definition}

\begin{remark}
\label{determined}
If $\tilde\Sigma'$, $\tilde h$ is good for $\Sigma'$, $h$, then
$\tilde h$ and $\tilde \Sigma'$ are completely determined by $h$ and
the values of $\tilde h$ on $C(\nabla^{\check h'}\times\{1\})$. Indeed,
$\tilde\Sigma'$ is a subdivision of $\tilde\Sigma$, and all one-dimensional
cones of $\tilde\Sigma'$ are either one-dimensional cones of $\Sigma'$ or
contained in $C(\nabla^{\check h'}\times\{1\})$. Thus $\tilde\Sigma'$ is
the unique such coherent subdivision on which the given values of $\tilde h$
extend to a convex piecewise linear function on $\tilde\Sigma'$. More
precisely, if we have a cone $\rho=C(\sigma)\times\{0\}+C(\tau\times\{1\})$
of $\tilde\Sigma$ with $\sigma$ a face of $\Delta^*$, $\tau$ a face of
$\nabla^{\check h'}$, the values of $\tilde h$ are given on $C(\sigma)\times
\{0\}$ and $C(\tau\times\{1\})$. Let
\[
\tilde\rho=Conv\left(\{(n,0,h(n))|n\in C(\sigma)\}\cup\{(n,l,\tilde h(n,l))|
(n,l)\in C(\tau\times\{1\})\}\right).
\]
Then $\tilde h$ can be defined by, for $(n,l)\in\rho$,
\[
\tilde h(n,l)=\inf\{r| (n,l,r)\in\tilde \rho\}.
\]
Of course, $\tilde h$ is strictly convex 
piecewise linear on a coherent subdivision
of $\rho$ which does not involve any one-dimensional rays not in 
$C(\sigma)\times\{0\}$ or $C(\tau\times\{1\})$, and this is the unique
such extension given the values on $C(\sigma)$ and $C(\tau\times\{1\})$.

One should note that an arbitrary choice of values for $\tilde h$ on
$C(\nabla^{\check h'}\times\{1\})$ does not necessarily yield a function
$\tilde h$ which is strictly convex on the whole fan $\tilde\Sigma'$,
and this can make it difficult, in general, to find good choices for
$\tilde h$.
\end{remark}
\bigskip

All of the above story can be repeated, interchanging the roles of 
$\Delta$ and $\nabla$, $h$ and $\check h$, so we get fans $\widetilde{\check\Sigma}$
and the notion of good subdivisions ${\widetilde{\check\Sigma'}}$, etc.

Now since the choice of $\tilde\Sigma'$, $\tilde h$, $\widetilde{\check\Sigma'}$,
$\widetilde{\check h}$ are not unique (even their existence hasn't been
addressed yet, see the end of this section), 
there is no reason for the corresponding dual intersection
complexes $(B_{\nabla}^{\check h},\P_{\nabla}^{\tilde\Sigma'})$
and $(B^h_{\Delta},\P^{\widetilde{\check\Sigma'}}_{\Delta})$ to be related
by a discrete Legendre transform. However, we will show for any good choice
of $\tilde\Sigma'$, $\tilde h$, we can find the correct 
$\widetilde{\check\Sigma'}$,
$\widetilde{\check h}$ which will yield dual intersection complexes related
by a discrete Legendre transform. This will be done essentially via a
Legendre transform in $N_{\RR}$.

Suppose we are given 
$\tilde\Sigma'$, $\tilde h$ good for $\Sigma',h$.
We have
\[
\tilde\nabla=\{(n,l)|n\in\nabla, l\ge h'(n)\}\subseteq N_{\RR}\oplus\RR.
\]
Set
\[
\tilde\nabla'=\{(n,l)|n\in\nabla^{\check h'}, l\ge\tilde h'(n,1)\}\subseteq
N_{\RR}\oplus\RR.
\]
Let $\widetilde{\check\Sigma'}$ be the normal fan to $\tilde\nabla+\tilde\nabla'$
in $M_{\RR}\oplus\RR$, and set
\[
\widetilde{\check h}(m,l)=-\inf\{\langle (m,l),(n,l')\rangle| (n,l')\in\tilde\nabla
+\tilde\nabla'\}
\]
to be the piecewise linear function on $\widetilde{\check\Sigma'}$ defined by
$\tilde\nabla+\tilde\nabla'$. Let $\widetilde{\check h'}=\widetilde{\check h}-
\widetilde{\check\varphi}$, so
\[
\widetilde{\check h'}(m,l)=-\inf\{\langle (m,l),(n,l')\rangle| (n,l')\in
\tilde\nabla'\}.
\]
Set
\[
\tilde\Delta'=\{(m,l)| m\in\Delta^{h'}, l\ge \widetilde{\check h'}(m,1)\}.
\]

We first wish to show that $\widetilde{\check\Sigma'}, \widetilde{\check h}$
is good for $\check\Sigma'$, $\check h$ and that this setup is 
dual, i.e. $\tilde\Sigma'$, $\tilde h$ can be recovered from
$\widetilde{\check\Sigma'}$, $\widetilde{\check h}$ in the same manner.
To this end,

\begin{lemma}
\label{redescribe}
\[
\tilde\nabla+\tilde\nabla'=\{(n,l)|\hbox{$n\in\nabla^h$ and
$l\ge \tilde h'(n,1)$}\}.
\]
\end{lemma}

Proof. First we need

\begin{claim}
Let $C(\sigma)\times\{0\}+C(\tau\times\{1\})$ be
a cone of $\tilde\Sigma'$, $\sigma\subseteq \partial\Delta^*$, $\tau
\subseteq\partial\nabla^{\check h'}$. Then
\[
(C(\sigma)\cap\nabla)+\tau=(C(\sigma)+\tau)\cap\nabla^{\check h}.
\]
\end{claim}

Proof. Inclusion of the first set in the second is obvious. Conversely, as in
the proof of Claim \ref{mink1}, it is enough to show that if $\rho$
is a proper face of $(C(\sigma)\cap\nabla)+\tau$ not contained in $\partial
(C(\sigma)+\tau)$, then $\rho\subseteq\partial\nabla^{\check h}$. Let
$\rho_i=\beta_i^*(\sigma)$, $\sigma_i=Conv\{0,\rho_i\}$.
By Claim \ref{mink1}, $C(\sigma)\cap\nabla=\sigma_1+\cdots+\sigma_r$.
Furthermore, if $\rho$ is a face as above, then there is some $i$
such that 
\[
\rho\subseteq \sigma_1+\cdots+\rho_i+\cdots+\sigma_r+\tau.
\]
Now $\sigma\subseteq\partial\Delta^*$, so there exists a face
$\omega\subseteq\partial\Delta$ dual to the minimal face of
$\Delta^*$ containing $\sigma$, and there exists a cone
$\tau'\in \check\Sigma'_{\check h'}$ corresponding to the minimal face
of $\nabla^{\check h'}$ containing $\tau$. By Lemma \ref{corrFh},
$\omega\cap\tau'\not=\emptyset$. In particular, there is a subcone
$\tau''\in\check\Sigma'$ of $\tau'$ with $\omega\cap\tau''\not=
\emptyset$. Let $m\in \omega\cap\tau''$. It follows from Claim \ref{mink1}
that $m=m_1+\cdots+m_r$ for $m_j\in\tau''\cap\Delta_j$,
$\langle m_i,\sigma_j\rangle=0$ for $i\not=j$, $\langle m_i,\rho_i\rangle
=-1$, and $\langle m_i,\tau\rangle =-\check h'(m_i)$. Thus
\[
\langle m_i,\rho\rangle =-1-\check h'(m_i)=-\check h(m_i).
\]
Thus we conclude that $\rho\subseteq\partial\nabla^{\check h}$.
\qed

Now define $\bar h:\nabla^h\rightarrow\RR$ by 
\[
\bar h(n)=\min\{l|(n,l)\in \tilde\nabla+\tilde\nabla'\}.
\]
Then 
\[
\bar h(n)=\min\{h'(n_1)+\tilde h'(n_2,1)| n_1\in \nabla,
n_2\in\nabla^{\check h'}, n_1+n_2=n\}.
\]
By convexity of $\tilde h'$, $\bar h(n)\ge \tilde h'(n,1)$.
If $n\in\nabla^{\check h'}$, then
\[
\bar h(n)\le h'(0)+\tilde h'(n,1)=\tilde h'(n,1),
\]
so $\bar h=\tilde h'(\cdot,1)$ on $\nabla^{\check h'}$. On the other
hand, if $n\in \nabla^{\check h}\setminus \nabla^{\check h'}$, then
$(n,1)$ is contained in a cone $C(\sigma)\times\{0\}+
C(\tau\times\{1\})$, and thus we can write, by the claim,
$n=n_1+n_2$ for $n_1\in C(\sigma)\cap\nabla$,
$n_2\in\tau$. Then since $\tilde h'$ is linear on this cone,
\[
\bar h(n)\le h'(n_1)+\tilde h'(n_2,1)=\tilde h'(n,1).
\]
Thus $\bar h=\tilde h'(\cdot, 1)$ on $\nabla^{\check h}$, as desired.
\qed

\begin{theorem}
\label{dualgood}
$\widetilde{\check\Sigma'}$ and $\widetilde{\check h}$ are good for 
$\check\Sigma'$,
$\check h$.
\end{theorem}

Proof. Because $\widetilde{\check\Sigma}$ is the normal fan to $\tilde\nabla$,
$\widetilde{\check\Sigma'}$ is a subdivision of $\widetilde{\check\Sigma}$. 
To check that it is good for $\check\Sigma'$, we need to consider normal
vectors to maximal lower faces of $\tilde\nabla+\tilde\nabla'$. Such a face
is a graph of $\tilde h'(\cdot,1)$ over its projection to $N_{\RR}$
by Lemma \ref{redescribe}. Let $\rho\in\tilde\Sigma'$ be a maximal 
cone, with $-\tilde h'$ given by $(m,l)\in M\oplus\ZZ$, and let
\[
\tilde\rho=\{(n,\tilde h'(n,1))| (n,1)\in\rho\cap (\nabla^{\check h}\times
\{1\})\},
\]
so that $\tilde\rho$ is contained in a maximal lower face of $\tilde
\Delta+\tilde\Delta'$. The primitive normal vector to this face is
then $(m,1)$, so we only need to check that $m\in\Delta^{h'}$. 
But by convexity of $\tilde h'$, 
\[
\langle m,n\rangle =\langle (m,l),(n,0)\rangle
\ge -\tilde h'(n,0)=-h'(n)
\]
for all $n\in N_{\RR}$, so $m\in \Delta^{h'}$.
\qed
\bigskip

The next result shows that by applying this procedure to 
$\widetilde{\check\Sigma'}$
and $\widetilde{\check h}$, we get back to $\tilde\Sigma'$, $\tilde h$.

\begin{theorem}
\label{secondtransform}
\[
\tilde h'(n,l)=-\inf\{\langle(m,l'),(n,l)\rangle| (m,l')\in\tilde\Delta'\}
\]
and
\[
\tilde h(n,l)=-\inf\{\langle(m,l'),(n,l)\rangle| 
(m,l')\in\tilde\Delta+\tilde\Delta'\}
\]
for any $(n,l)$ in the support of the fan $\tilde\Sigma'$.
\end{theorem}

Proof.
By the dual of Theorem \ref{dualgood}, $\tilde\Delta+\tilde\Delta'$
induces a good subdivision $\tilde\Sigma_1'$, $\tilde h_1$ for $\Sigma',
h$. We just need to show $\tilde h_1=\tilde h$, but by Remark \ref{determined},
it is enough to check this on $C(\nabla^{\check h'}\times\{1\})$.
But $\tilde h=\tilde h'$ on $C(\nabla^{\check h'}\times\{1\})$, so it
is enough to check the first equality on $\nabla^{\check h'}\times\{1\}$.

Let $n\in\nabla^{\check h'}$, $(m,l)\in\tilde\Delta'$. Then
\begin{eqnarray*}
\langle (m,l),(n,1)\rangle &=& \langle m,n\rangle +l\\
&\ge& \langle m,n\rangle +\widetilde{\check h'}(m,1)\\
&=&\langle m,n\rangle -\inf \{\langle (m,1), (n',l')\rangle | (n',l')\in
\tilde\nabla'\}\\
&\ge&\langle m,n\rangle -\langle (m,1),(n,\tilde h'(n,1))\rangle\\
&=&-\tilde h'(n,1).
\end{eqnarray*}
Thus
\[
\tilde h'(n,1)\ge -\inf\{\langle (m,l),(n,1)\rangle | (m,l)\in
\tilde\nabla'\}.
\]
Conversely, suppose $n\in\nabla^{\check h'}$ is such that $(n,1)$ is in
a maximal cone of $\tilde\Sigma'$ contained in $C(\nabla^{\check h'}\times
\{1\})$. Then on this cone $-\tilde h'$ is given by some
$(m,l)\in M_{\RR}\oplus \RR$, so $l=-\tilde h'(n,1)-\langle m,n\rangle$, and
\begin{eqnarray*}
\widetilde{\check h'}(m,1)&=&-\inf\{\langle (m,1), (n',l')\rangle|
(n',l')\in\tilde\nabla'\}\\
&=&-\inf\{\langle (m,1),(n',\tilde h'(n',1))\rangle| n'\in\nabla^{\check h'}\}\\
&=&-\inf\{\langle m,n'\rangle+\tilde h'(n',1)|n'\in\nabla^{\check h'}\}.
\end{eqnarray*}
By convexity of $\tilde h'$, this infinum occurs when $n'=n$, so
$\widetilde{\check h'}(m,1)=l$.
On the other hand, as observed in the proof of Theorem \ref{dualgood},
$m\in\Delta^{h'}$. Thus $(m,l)\in\tilde\Delta'$ and
\[
-\inf\{\langle(m',l'),(n,1)\rangle|(m',l')\in\tilde\nabla'\}
\ge -\langle (m,l),(n,1)\rangle=\tilde h'(n,1).
\]
Thus equality holds. \qed

\bigskip
Thus we have a way of passing between good data $\tilde\Sigma'$,
$\tilde h$ and good data $\widetilde{\check\Sigma'}$, $\widetilde{\check h}$,
and this procedure is reversible. To complete the story, we need
a method of passing between relevant cones of $\tilde\Sigma'$ and relevant
cones of $\widetilde{\check\Sigma'}$.

Let 
\[
\rho=C(\sigma)\times\{0\}+C(\tau\times\{1\})
\]
be a relevant cone of $\tilde\Sigma'$. Then $\lambda(\rho):=\mbeta(\sigma)+\tau$
is contained in a face of $\nabla^h$. Let
\begin{eqnarray*}
\tilde\lambda(\rho)&:=&\{(n,h'(n))| n\in\mbeta(\sigma)\}+
\{(n,\tilde h'(n,1))| n\in\tau\}\\
&=&\{(n,\tilde h'(n,1))| n\in\lambda(\rho)\},
\end{eqnarray*}
the last equality following from $\tilde h'$ being linear on $\rho$.
By Lemma \ref{redescribe}, $\tilde\lambda(\rho)\subseteq\tilde\nabla
+\tilde\nabla'$.
Let
\[
\alpha(\rho)=N_{\tilde\nabla+\tilde\nabla'}(\tilde\lambda(\rho)).
\]
Then $\alpha(\rho)$ is a cone in $\widetilde{\check\Sigma'}$. We can
similarly, for a relevant cone $\check\rho$ of $\widetilde{\check\Sigma'}$,
define $\check\lambda(\check\rho)$, $\widetilde{\check\lambda}(\check\rho)$, 
and a cone $\check\alpha(\check\rho)$ of $\tilde\Sigma'$.

\begin{proposition} 
\label{relevance}
If $\rho$ is a relevant cone of $\tilde\Sigma'$,
then $\widetilde\lambda(\rho)$ is a face of $\tilde\nabla+\tilde\nabla'$
and $\alpha(\rho)$ is a relevant cone of $\widetilde{\check\Sigma'}$.
\end{proposition}

Proof. First observe that $\alpha(\rho)$ is 
a cone of the second type of Observation \ref{conetypes}. To see this,
just observe that $\tilde h'(\cdot,1)$ is affine linear on
$\mbeta(\sigma)+\tau$ and $\mbeta(\sigma)+\tau$ is contained in
$\partial\nabla^{\check h}$, so by Lemma \ref{redescribe}, 
$\tilde\lambda(\rho)$ is contained in a face
of $\tilde\nabla+\tilde\nabla'$ projecting injectively into
$\partial\nabla^{\check h}$. Then $\alpha(\rho)$
is clearly of the second type.

Note now that the faces of $\tilde\nabla+\tilde\nabla'$ whose normal
cones are of the second type are in one-to-one correspondence via 
projection $\check p:N_{\RR}\oplus\RR\rightarrow N_{\RR}$ 
with the collection of subsets of $\partial
\nabla^{\check h}$ given by
\[
\{\check p((\omega\times\{1\})\cap\delta)|\hbox{$\omega$ is
a face of $\nabla^{\check h}$ and $\delta$ is a normal cone
to $\tilde\Delta'$}\}.
\]
Indeed, the function $\tilde h'$ is linear on normal cones of $\tilde
\Delta'$ by Theorem \ref{secondtransform}, and the statement
then follows from Lemma \ref{redescribe}. So
to show $\widetilde\lambda(\rho)$ is a face of $\tilde\nabla+\tilde\nabla'$
we need to show $\lambda(\rho)=\check p((\omega\times\{1\})\cap
\delta)$ for some face $\omega$ of $\nabla^{\check h}$ and $\delta$
a normal cone of $\tilde\Delta'$. But as $\rho$ is a normal cone to
a face of $\tilde\Delta+\tilde\Delta'$, it is necessarily of the
form $\rho'\cap\delta$, where $\rho'$ is a normal cone to
$\tilde\Delta$, i.e. $\rho'\in\tilde\Sigma$, and $\delta$ is a normal
cone to $\tilde\Delta'$. Then $\rho'=C(\sigma')\times\{0\}+
C(\tau'\times\{1\})$ with $\sigma'+\tau'$ a face of 
$\Delta^*+\nabla^{\check h'}$, and $\mbeta(\sigma')+\tau'=:\omega$ is a face of
$\nabla+\nabla^{\check h'}=\nabla^{\check h}$. Using the notation of the proof
of Theorem \ref{toricdegen2}, we can take $V_{\rho}=V_{\rho'}$, and then
by Claim \ref{basicclaim},
\begin{eqnarray*}
\lambda(\rho)\times\{1\}&=&(\mbeta(\sigma)+\tau)\times\{1\}\\
&=&V_{\rho}^{\perp}\cap\rho\cap(N_{\RR}\times\{1\})\\
&=&V_{\rho}^{\perp}\cap\rho'\cap\delta\cap(N_{\RR}\times\{1\})\\
&=&(\omega\times\{1\})\cap\delta,
\end{eqnarray*}
as desired.

Next write $\check\rho=\alpha(\rho)$,
\[
\check\rho=C(\check\sigma)\times\{0\}+C(\check\tau\times\{1\}).
\]
To show $\check\rho$ is relevant, we show $\beta(\check\sigma)\not=\emptyset$.
Now $\rho$ is of course contained in $\rho'\in\tilde\Sigma$, 
so $\sigma\subseteq\sigma'$ and
$\tau\subseteq\tau'$. Changing the meaning of $\omega$,
for notational convenience,
write $\omega^*:=\sigma'$, so we can write $\mbeta(\sigma')=\momega$,
and write $\omega,\momega^*$ as in \S 2.
Then $\mbeta(\sigma)\subseteq\momega$. 

Now $C(\check\sigma)=N_{\nabla^{\check h}}(\mbeta(\sigma)+\tau),$
and as $\nabla^{\check h}=\nabla+\nabla^{\check h'}$,
\[
C(\check\sigma)=N_{\nabla}(\mbeta(\sigma))\cap N_{\nabla^{\check h'}}
(\tau),
\]
$C(\momega^*)=N_{\nabla}(\momega)\subseteq N_{\nabla}(\mbeta(\sigma))$,
and $N_{\nabla^{\check h'}}(\tau')\subseteq N_{\nabla^{\check h'}}(\tau)$.
On the other hand, $\omega=\beta(\momega^*)$ is the dual face of $\Delta$
to $\omega^*$ by Lemma \ref{corrbetweenfaces}, and $\omega^*+\tau'$
is a face of $\Delta^*+\nabla^{\check h'}$ by Proposition \ref{sigmatwiddle}.
Thus
\begin{eqnarray*}
\emptyset\not=N_{\Delta^*+\nabla^{\check h'}}(\omega^*+\tau')&=&
N_{\Delta^*}(\omega^*)\cap N_{\nabla^{\check h'}}(\tau')\\
&=&C(\omega)\cap N_{\nabla^{\check h'}}(\tau').
\end{eqnarray*}
Since $C(\omega)\subseteq C(\momega^*)$ and $\check\varphi_i$ is non-zero
at each non-zero point of $C(\omega)$, it follows that for each
$i$, there exists a point in $C(\momega^*)\cap N_{\nabla^{\check h'}}(\tau')$
on which $\check\varphi_i$ is non-zero. Thus in particular, as
$\momega^*\cap N_{\nabla^{\check h'}}(\tau')$ has integral vertices,
$\mbeta(\momega^*\cap N_{\nabla^{\check h'}}(\tau'))\not=\emptyset$.
But $\momega^*\cap N_{\nabla^{\check h'}}(\tau')\subseteq\check\sigma$,
so $\beta(\check\sigma)\not=\emptyset$. \qed

\begin{proposition} $\alpha$ and $\check\alpha$ give a one-to-one
order reversing correspondence between relevant cones of
$\tilde\Sigma'$ and relevant cones of $\widetilde{\check\Sigma'}$.
\end{proposition}

Proof. We need to show that $\rho=\check\alpha(\alpha(\rho))$. First,
we show 
\[
\rho\subseteq\check\alpha(\alpha(\rho))=N_{\tilde\Delta+\tilde
\Delta'}(\widetilde{\check\lambda}(\alpha(\rho))).
\] 
Write $\rho=C(\sigma)\times\{0\}+C(\tau\times\{1\})$, $\check
\rho=\alpha(\rho)$. Let
$n\in\beta_i^*(\sigma)\subseteq\nabla_i$. Then $(n,h'(n))$ coincides with
the function $-\widetilde{\check\varphi}_i$ on $\check\rho
=C(\check\sigma)\times\{0\}+C(\check\tau\times\{1\})$.
Now if $m_1\in\beta(\check\sigma)$, $-\widetilde{\check\varphi}_i(m_1,0)=-1$,
and if $m_2\in\check\tau$, $-\widetilde{\check\varphi}_i(m_2,1)=0$. Thus
\begin{eqnarray*}
\langle (n,0),(m_1,\check h'(m_1))+(m_2,\widetilde{\check h'}(m_2,1))\rangle
&=&\langle (n,h'(n)),(m_1,0)+(m_2,1)\rangle-h'(n)\\
&=&-1-h'(n),
\end{eqnarray*}
so $(n,0)$ is constant on $\widetilde{\check\lambda}(\alpha(\rho))$.
If on the other hand $(m_1,l_1)\in\tilde\Delta$, $(m_2,l_2)\in
\tilde\Delta'$, then
\begin{eqnarray*}
\langle (n,0),(m_1,l_1)+(m_2,l_2)\rangle
&=&\langle (n,h'(n)),(m_1,0)+(m_2,1)\rangle -h'(n)\\
&\ge&-\widetilde{\check\varphi}_i(m_1,0)-\widetilde{\check\varphi}_i(m_2,1)- h'(n)\\
&\ge&-1-h'(n)
\end{eqnarray*}
Thus $(n,0)\in N_{\tilde\Delta+\tilde\Delta'}(\widetilde{\check\lambda}
(\alpha(\rho)))$.
Just as in the proof of Lemma \ref{corrbetweenfaces},
$\sigma=Conv\{\beta_1^*(\sigma),\ldots,\beta_r^*(\sigma)\}$
so $C(\sigma)\times\{0\}\subseteq\check\alpha(\alpha(\rho))$.

Next let $n\in\tau$. 
Consider $(m_1,l_1)\in\tilde\Delta$, $(m_2,l_2)\in\tilde\Delta'$.
Now by Theorem \ref{secondtransform},
\[
\langle (n,1),(m_1,l_1)+(m_2,l_2)\rangle \ge -\tilde h'(n,1).
\]
On the other hand, by the definition of $\widetilde{\check h}$,
$(n,\tilde h'(n,1))$ defines the function
$-\tilde{\check h}'$ on $\check\rho$, so if $m_1\in\beta(\check\sigma)$, 
$m_2\in\check\tau$, 
\[
\langle (n,\tilde h'(n,1)),(m_1,0)+(m_2,1)\rangle=-\tilde {\check h}'(m_1+m_2,1)
\]
or
\[
\langle n,m_1+m_2\rangle+\widetilde{\check h'}(m_1+m_2,1)=-\tilde h'(n,1)
\]
or
\[
\langle (n,1),(m_1,\check h'(m_1))+(m_2,\widetilde{\check h'}(m_2,1))\rangle 
=-\tilde h'(n,1).
\]
Thus $(n,1)$ takes the constant value $-\tilde h'(n,1)$
on $\widetilde{\check\lambda}(\alpha(\rho))$, so  $C(\tau\times\{1\})\subseteq
\check\alpha(\alpha(\rho))$. Thus
$\rho\subseteq\check\alpha(\alpha(\rho))$.

To show equality, we compute dimensions. Let $n=\dim M_{\RR}$.
Then if $\dim\rho=k$, $\rho$ corresponds to an $n+1-k$ dimensional
stratum of $X(\tilde\Sigma')$, and hence $\lambda(\rho)$ is 
an $(n-r)-(n+1-k)$ dimensional cell of $\P_{\nabla}^{\tilde\Sigma'}$,
and $\widetilde\lambda(\rho)$ is the same dimension. By Proposition \ref{relevance},
$\widetilde\lambda(\rho)$ is a face of $\tilde\nabla+\tilde\nabla'$, so
$\dim\alpha(\rho)=(n+1)-[(n-r)-(n+1-k)]=n-k+r+2$. Then 
$\dim\check\alpha(\alpha(\rho))=n-(n-k+r+2)+r+2=k$, so
$\rho=\check\alpha(\alpha(\rho))$. \qed
\bigskip

The functions $\tilde h$ and $\widetilde{\check h}$ give ample divisors on
$X(\tilde\Sigma')$ and $X(\widetilde{\check\Sigma'})$ respectively, and these
restrict to polarizations on $\X'\rightarrow \shS$ and $\check\X'\rightarrow
\shS$. These then give rise to multi-valued piecewise linear functions
$\tilde h_{\nabla}$ and $\widetilde{\check h}_{\Delta}$ on $B_{\nabla}^{\check h}$
and $B_{\Delta}^h$ respectively. Finally, we have

\begin{theorem}
$(B_{\Delta},\P_{\Delta}^{\widetilde{\check\Sigma'}},
\widetilde{\check h}_{\Delta})$ is the discrete Legendre transform of 
$(B_{\nabla}^{\check h},\P_{\nabla}^{\tilde\Sigma'},\tilde h_{\nabla})$.
\end{theorem}

Proof. By symmetry of the situation, it is enough to check that
given a relevant cone $\rho$ of $\tilde\Sigma'$, the piecewise
linear function on the quotient fan $\tilde\Sigma'(\rho)$ induced
by $\tilde h$ has as 
Newton polytope $\check\lambda(\check\alpha(\rho))$. But
$\check\lambda(\check\alpha(\rho))$ can be identified via projection
with $\widetilde{\check\lambda}(\check\alpha(\rho))$, and by
Theorem \ref{secondtransform}, the latter is canonically 
identified with this Newton polytope. \qed
\bigskip

Finally we discuss existence of pairs $\tilde h, \tilde\Sigma'$ good for
$h,\Sigma'$. Suppose that $h$ and $\check h$ are given. Now I do
not expect that in general one can find a $\tilde h,\tilde\Sigma'$
good for $h, \Sigma'$, but I do conjecture that one can find
a pair $\tilde h,\tilde\Sigma'$ good for $nh,\Sigma'$, for $n$ a 
sufficiently large positive integer. This affects the affine
structure on the intersection complex merely by rescaling by the factor
$n$, so this is not a serious restriction. It arises due to integrality
issues, as we will see below. However, this still seems to be a 
difficult combinatorial problem which has so far resisted attempts
to solve. Thus we will settle for a weaker statement:

\begin{theorem}
There exists positive integers $n_0$ and $m_0$ such that for all $n\ge n_0$,
there exists $\tilde h,\tilde\Sigma'$ good for $m_0h+n\varphi, \Sigma$.
\end{theorem}

Proof. Begin by choosing the values of $\tilde h$ on $C(\nabla^{\check h'}\times
\{1\})$. This can be done essentially arbitrarily; for convenience we
take $\tilde h\equiv 0$ on this cone. In addition, we wish to take
$\tilde h=h$ on $N_{\RR}\times\{0\}$. Using Remark \ref{determined}, 
each cone $\rho=C(\sigma)\times\{0\}+C(\tau\times\{1\})$
of $\tilde\Sigma$ can be subdivided in such a way that $\tilde h$ extends
to a strictly convex piecewise linear function on this subdivision of
$\rho$. These subdivisions are compatible: if $\rho\subseteq\rho'$ and
we subdivide $\rho$ and $\rho'$ as in Remark \ref{determined}, then
the subdivision of $\rho'$ restricts to the subdivision of $\rho$.
Hence we obtain a subdivision $\tilde\Sigma'$ of $\tilde\Sigma$
and a function $\tilde h$ on $\tilde\Sigma'$. It is not
necessarily integral, but does have rational slopes, so there exists
$m_0\in\ZZ$ such that $m_0\tilde h$ is integral. Now strict 
convexity can be tested across codimension one cones of the fan
$\tilde\Sigma'$, and by construction the only cones which might
cause problems are codimension one cones of $\tilde\Sigma'$ contained
in codimension one cones of $\tilde\Sigma$. But $\tilde\varphi$ is
strictly convex on $\tilde\Sigma$, so there exists an $n_0$ such that
$m_0\tilde h+n\tilde\varphi$ is convex on $\tilde\Sigma'$ for
all $n\ge n_0-1$ and strictly convex for $n\ge n_0$. 
Thus $m_0\tilde h+n\tilde\varphi$, $\tilde \Sigma'$
is good for $m_0h+n\varphi$, $\Sigma'$.
\qed
\bigskip

We note that as $h$ defines an ample divisor on $\PP_{\Delta^h}$, $m_0h+n\varphi$
is also ample. Thus from the point of view of applying Theorem \ref{torusbundle}
to understand the topology of the SYZ fibration, this result is sufficient in the light of
Theorem \ref{simplicity}.

\section{An example, and connection to Kovalev's example}

We will consider Chad Schoen's example \cite{Schoen} of a Calabi-Yau threefold
obtained as a fibred product of two rational elliptic surfaces. Let
$f_1:Y_1\rightarrow\PP^1$ and $f_2:Y_2\rightarrow\PP^1$ be two rational
elliptic surfaces with section, and suppose there is no point
$x\in\PP^1$ such that both $f_1^{-1}(x)$ and $f_2^{-1}(x)$ 
are singular curves. Then $X=Y_1\times_{\PP^1} Y_2$ is a non-singular
Calabi-Yau threefold with $\chi(X)=0$ and $h^{1,1}(X)=h^{1,2}(X)=19$. 

To study mirror symmetry for the Schoen fibred product, 
we proceed as in \cite{HSS} by representing
it as a complete intersection in $\PP^1\times\PP^2\times\PP^2$
of hypersurfaces of tridegree $(1,3,0)$ and $(1,0,3)$. In terms of
toric data, set
\begin{eqnarray*}
\bar\Delta_1&=&[-1,1]\subseteq\RR\\
\bar\Delta_2=\bar\Delta_3&=&Conv\{(-1,-1),(2,-1),(-1,2)\}\subseteq\RR^2\\
\bar\Delta_1^+&=&[0,1]\subseteq\RR\\
\bar\Delta_1^-&=&[-1,0]\subseteq\RR
\end{eqnarray*}
and then
\begin{eqnarray*}
\Delta&=&\bar\Delta_1\times\bar\Delta_2\times\bar\Delta_3
\subseteq \RR\times\RR^2\times\RR^2=M_{\RR}\\
\Delta_1&=&\bar\Delta_1^+\times\bar\Delta_2\times\{0\}\\
\Delta_2&=&\bar\Delta_1^-\times\{0\}\times\bar\Delta_3.
\end{eqnarray*}
Label the vertices of $\Delta_1$ by
\begin{eqnarray*}
P_0^0=(0,-1,-1,0,0),& P_1^0=(0,2,-1,0,0),& P_2^0=(0,-1,2,0,0)\\
P_0^+=(1,-1,-1,0,0),& P_1^+=(1,2,-1,0,0),& P_2^+=(1,-1,2,0,0)
\end{eqnarray*}
and the vertices of $\Delta_2$ by 
\begin{eqnarray*}
Q_0^0=(0,0,0,-1,-1),& Q_1^0=(0,0,0,2,-1),& Q_2^0=(0,0,0,-1,2)\\
Q_0^-=(-1,0,0,-1,-1),& Q_1^-=(-1,0,0,2,-1),& Q_2^-=(-1,0,0,-1,2)
\end{eqnarray*}
so that 
\[
\nabla^*=Conv\{P_i^0,P_i^+,Q_i^0,Q_i^-|i=0,1,2\}.
\]
Setting $\bar\Delta_1^*=[-1,1]\subseteq\RR$ and
\[
\bar\Delta_2^*=\bar\Delta_3^*=Conv\{(-1,-1),(1,0),(0,1)\}\subseteq\RR^2,
\]
then one sees easily that 
\begin{eqnarray*}
\Delta^*&=&Conv\{\bar\Delta_1^*\times\{0\}\times\{0\},
\{0\}\times\bar\Delta_2^*\times\{0\},\{0\}\times\{0\}\times\bar\Delta_3^*\}\\
&\subseteq&\RR\times\RR^2\times\RR^2=N_{\RR}.
\end{eqnarray*}
Let $R^+=(1,0,0,0,0)$ and $R^-=(-1,0,0,0,0)$,
\begin{eqnarray*}
S_0=(0,-1,-1,0,0),& S_1=(0,1,0,0,0),& S_2=(0,0,1,0,0)\\
T_0=(0,0,0,-1,-1),& T_1=(0,0,0,1,0),& T_2=(0,0,0,0,1)
\end{eqnarray*}
be points in $N_{\RR}$. Then
\begin{eqnarray*}
\nabla_1&=&Conv\{0,R^-,S_0,S_1,S_2\}\\
\nabla_2&=&Conv\{0,R^+,T_0,T_1,T_2\}
\end{eqnarray*}
and one can compute that
\[
\nabla=\nabla_1+\nabla_2=Conv\{R^-+T_j, R^++S_i, S_i+T_j|i,j=0,1,2\}.
\]
Finally, we take $h=\varphi$ (corresponding to the anti-canonical
polarization on $\PP_{\Delta}=\PP^1\times\PP^2\times\PP^2$) and choose
$\check h$ to give an ample divisor on a partial resolution of $\PP_{\nabla}$.
We will discuss the effect of this latter choice shortly. We can now
describe $B_{\Delta^h}=B_{\Delta}$ as a polyhedral subcomplex of
$\Delta$. There is one three-dimensional cell of $B_{\Delta}$
corresponding to each vertex
of $\nabla$, and these are of the form (taking indices modulo 3)
\begin{eqnarray*}
\sigma_j^+&=&Conv\{P_0^+,P_1^+,P_2^+\}+Conv\{Q_{j+1}^0,Q_{j+2}^0\}\\
\sigma_i^-&=&Conv\{P_{i+1}^0,P_{i+2}^0\}+Conv\{Q_0^-,Q_1^-,Q_2^-\}\\
\sigma_{ij}&=&Conv\{P_{i+1}^0,P_{i+2}^0,P_{i+1}^+,P_{i+2}^+\}
+Conv\{Q_{j+1}^0,Q_{j+2}^0,Q_{j+1}^-,Q_{j+2}^-\}
\end{eqnarray*}
corresponding to $R^-+T_j$, $S_i+R^+$ and $S_i+T_j$ respectively.
Then 
\begin{eqnarray*}
\bigcup_{j=0,1,2} \sigma_j^+&=&\{1\}\times\bar\Delta_2\times\partial
\bar\Delta_3\\
\bigcup_{i=0,1,2} \sigma_i^-&=&\{-1\}\times\partial\bar\Delta_2\times
\bar\Delta_3
\end{eqnarray*}
are both solid tori,
and
\[
\bigcup_{i,j=0,1,2}\sigma_{ij}=\bar\Delta_1\times\partial\bar\Delta_2
\times\partial\bar\Delta_3.
\]
These three solids glue together to form an $S^3$ (this is just
a variant of the usual Heegard splitting of the sphere).

We will describe the discriminant locus $\Gamma_{\Delta}$.
In previous sections, we took the discriminant locus $\Gamma$
to be pretty large, including all codimension two simplices of
$\Bar(\P_{\Delta}^{\widetilde{\check\Sigma'}})$ not containing
vertices of $\P_{\Delta}^{\widetilde{\check\Sigma'}}$ or intersecting
interiors of maximal cells. However, by Proposition 1.27 of \cite{AffI},
the affine structure can be extended across any simplex of
$\Gamma$ for which the monodromy about this simplex is trivial.
Thus, by applying Proposition \ref{monodromy2}, we can identify the
minimal discriminant locus.
This will depend on the choice of $\check h$,
which we do not wish to specify precisely. However, we can specify
the decomposition of $\nabla^*$ it induces, or less
specifically, the decomposition of certain faces of $\nabla^*$,
which will be enough for us.

Look first at $\sigma_i^+$. We can take, for example, the polyhedral
decompositions $\P_1$ and $\P_2$ of $\{1\}\times\bar\Delta_2\times\{0\}$
and $Conv\{Q_{j+1}^0,Q_{j+2}^0\}$ depicted here:
\begin{center}
\includegraphics{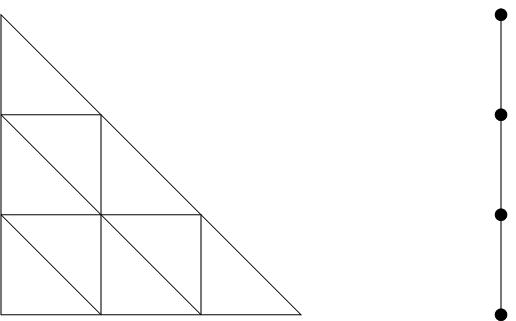}
\end{center}
Then we obtain a polyhedral decomposition of the face of 
$\nabla^*$ dual to $R^-+T_j$ (this face is $Conv\{P_0^+,P_1^+,P_2^+,
Q_{j+1}^0,Q_{j+2}^0\}$) given by
\[
\P_{R^-+T_j}=\{Conv(\tau_1\cup\tau_2)|\tau_1\in\P_1,\tau_2\in\P_2\}
\]
and then a polyhedral decomposition of $\sigma_j^+$ given by
\[
\P_{\sigma_j^+}=\{\tau_1+\tau_2|\tau_1\in\P_1,\tau_2\in\P_2\}.
\]
The polytopes of $\P_{\sigma_j^+}$ contained in $\partial\sigma_j^+$
are as depicted:
\begin{center}
\includegraphics{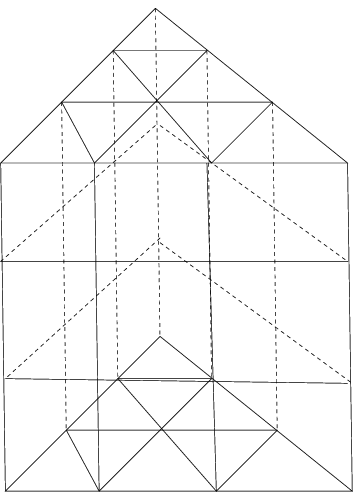}
\end{center}
As there is no monodromy around any loop contained in $Int(\sigma_j^+)$,
it is enough to focus on the boundary. By applying Proposition \ref{monodromy2},
one finds that the only bit of the discriminant locus that survives is
as depicted:
\begin{center}
\includegraphics{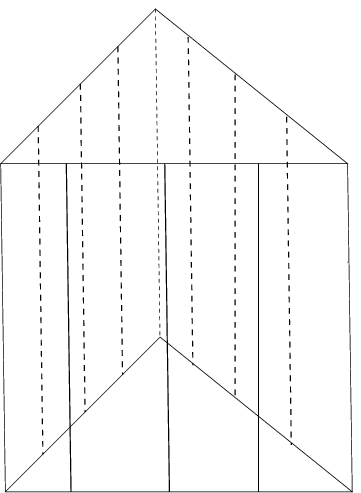}
\end{center}
We have a similar picture for $\sigma_i^-$. 

Finally let, for $i,k,l=0,1,2$,
\begin{eqnarray*}
{}^+\tau_i^k&=&Conv\bigg\{
{k\over 3} P^+_{i+1}+{3-k\over 3} P_{i+2}^+, {k+1\over 3} P_{i+1}^++
{2-k\over 3} P_{i+2}^+,\\
&&{k\over 3} P^0_{i+1}+{3-k\over 3} P_{i+2}^0, {k+1\over 3} P_{i+1}^0+
{2-k\over 3} P_{i+2}^0\bigg\},\\
{}^-\tau_i^k&=&Conv\bigg\{
{k\over 3} Q^-_{i+1}+{3-k\over 3} Q_{i+2}^-, {k+1\over 3} Q_{i+1}^-+
{2-k\over 3} Q_{i+2}^-,\\
&&{k\over 3} Q^0_{i+1}+{3-k\over 3} Q_{i+2}^0, {k+1\over 3} Q_{i+1}^0+
{2-k\over 3} Q_{i+2}^0\bigg\},\\
{}^+\rho_i^l&=&Conv\bigg\{{l\over 3} Q_{i+1}^0+{3-l\over 3} Q_{i+2}^0,
{l+1\over 3} Q^0_{i+1}+{2-l\over 3}Q^0_{i+2}\bigg\},\\
{}^-\rho_i^l&=&Conv\bigg\{{l\over 3} P_{i+1}^0+{3-l\over 3} P_{i+2}^0,
{l+1\over 3} P^0_{i+1}+{2-l\over 3}P^0_{i+2}\bigg\}.
\end{eqnarray*}
Then for a fixed $i,j$,
\[
\{Conv\{{}^+\tau_i^k,{}^+\rho_j^l\}| 0\le k,l\le 2\}
\cup
\{Conv\{{}^-\rho_i^k,{}^-\tau_j^l\}| 0\le k,l\le 2\}
\]
form the maximal cells of a polyhedral decomposition of the face
of $\nabla^*$ dual to the vertex $S_i+T_j$ of $\nabla$. Then this induces
a decomposition of $\sigma_{ij}$ with maximal cells of the
form 
\[
\{{}^+\tau_i^k+{}^+\rho_j^l| 0\le k,l\le 2\}
\cup
\{{}^-\rho_i^k+{}^-\tau_j^l| 0\le k,l\le 2\}.
\]
On the boundary of $\sigma_{ij}$ this looks like
\begin{center}
\input{sigmaij.pstex_t}
\end{center}
and the discriminant locus on the boundary of $\sigma_{ij}$ looks like
\begin{center}
\input{sigmaijdisc.pstex_t}
\end{center}

Putting this all together, we can describe the discriminant locus
as consisting of $24$ circles, occurring in two sets of 12. As a link, this
can be viewed as taking the Hopf link and replacing each circle of the Hopf
link with 12 parallel circles.

Now so far we have only followed through the combinatorics, and gotten
a specific dual intersection complex. The main point of this section is
to observe that this affine manifold, and a more general version of this
manifold,
can be obtained abstractly using surgeries of affine manifolds, inspired
by Kovalev's construction \cite{Kov}. The starting point is to find a structure of
integral affine manifold with singularities with boundary on a closed
2-disk $D$, so that the boundary is an affine line of integral length
and such that the affine structure on $D$ has 12 singular points in the 
interior with monodromy conjugate to $\begin{pmatrix}1&1\\ 0&1\end{pmatrix}$.
We leave as an exercise a construction of an example of such; for a hint,
either study carefully the construction above or consult
\cite{Sym}. The boundary of $D$ is an affine manifold isomorphic to
$\RR/n\ZZ$ for some integer $n$. Take $B_1=B_2=D\times(\RR/n\ZZ)$,
so $\partial B_1=\partial B_2=(\RR/n\ZZ)\times (\RR/n\ZZ)$, and glue
$\partial B_1$ and $\partial B_2$ along their boundaries using 
an involution $\iota$ of $(\RR/n\ZZ)\times (\RR/n\ZZ)$ given by $\iota(a,b)=
(b,a)$. After gluing, we get $B$, and because we are gluing $B_1$ and
$B_2$ across an affine plane, we can define an affine structure
on $B$ restricting to the given affine structures on $B_1$ and $B_2$.
With the correct choice of $D$, one can show one obtains the same
affine manifold which we constructed as a dual intersection complex 
(where incidentally $n=9$).

Now the point is that even if we didn't know about the Batyrev-Borisov
construction, we can easily produce a degeneration with $B$ as
dual intersection complex. Begin by finding a toric degeneration
of rational elliptic surfaces (well, we haven't given a definition
of such a degeneration, but use your imagination with Definition
\ref{toric degen} as a guide); in particular, we have a diagram
\[
\xymatrix@C=30pt
{\X\ar[d]^g\ar[rd]^f&\\
\PP^1\times\shS\ar[r]&\shS\\
}
\]

A typical example of such a degeneration is obtained in $\PP^1\times\PP^2
\times\shS$ with equation $tf+x_0y_0y_1y_2=0$ where $f$ is a general
polynomial of bidegree $(1,3)$ and $(x_0,x_1)$, $(y_0,y_1,y_2)$ are
coordinates on $\PP^1$ and $\PP^2$ respectively. This is actually
not quite general enough (we will obtain a disk $D$ below with,
instead of 12 singular points, only
three points of multiplicity four), but this should give enough
of a hint. As in this example, we will assume $g:\X_0\rightarrow\PP^1$
has as general fibre an $I_n$ fibre, and in particular $g^{-1}(\infty)$
is such a fibre, but $g^{-1}(0)$ can be a union of rational surfaces.
If we remove $g^{-1}(\{\infty\}\times\shS)$ from $\X$, we obtain a family
of open surfaces, and the dual intersection complex of this family makes
sense as a manifold with boundary: the fans of the open components of
$\X_0$ are not complete, and the corresponding vertices are in the
boundary of the dual intersection complex, which will be a disk $D$
with an affine structure with singularities. The boundary of $D$ is
$\RR/n\ZZ$.

To construct a degeneration of Schoen's fibred products, blow up the
point $\{\infty\}\times\{0\}$ in $\PP^1\times\shS$ and 
the curve $g^{-1}(\{\infty\}\times\{0\})$ in $\X$ to get 
\[
\xymatrix@C=30pt
{\tilde\X\ar[d]^g\ar[rd]^f&\\
\Y\ar[r]&\shS\\
}
\]
Now take another similar family $\X'\rightarrow\shS$, in which the role
of $0$ and $\infty$ in $\PP^1$ is reversed. We blow up $\{0\}\times\{0\}$
in $\PP^1\times\shS$ this time to get $\Y'$, and $(g')^{-1}(\{0\}\times
\{0\})$ in $\X'$ to give
\[
\xymatrix@C=30pt
{\tilde\X'\ar[d]^g\ar[rd]^f&\\
\Y'\ar[r]&\shS\\
}
\]
We can identify $\Y$ and $\Y'$ so that the exceptional curve on $\Y$ is
identified with the proper transform of $\PP^1\times\{0\}$ on $\Y'$,
and vice versa. We then take $\Z=\tilde\X\times_{\Y}\tilde\X'
\rightarrow\shS$ as a toric degeneration of Schoen's fibred product.
It is easy to see the dual intersection complex is constructed
by gluing two copies of
$D\times S^1$ with the twist as described above.

The decomposition of $B$ into $B_1$ and $B_2$ actually corresponds,
in some sense, to a simpler degeneration of the Schoen threefolds 
into a union of two threefolds, essentially as considered in
\cite{Kov}, on which one can explicitly find special Lagrangian fibrations.
In our algebro-geometric picture, we degenerate these two threefolds
further, to the threefolds sitting over the two curves in $\Y$ sitting
over $0\in\shS$. As a last observation, we note that it is easy to see
that dualizing the torus fibration obtained from $B$ doesn't change the
topology, using the fact that two-dimensional torus bundles are self-dual.
Hence we see an explicit prediction of the Strominger-Yau-Zaslow approach
that the Schoen fibred product is self-mirror, something which may not
have been self-evident from the Batyrev-Borisov picture.

\end{document}

%% file: sigmaij.pstex_t
\begin{picture}(0,0)%
\includegraphics{sigmaij.pstex}%
\end{picture}%
\setlength{\unitlength}{1579sp}%
\begingroup\makeatletter\ifx\SetFigFont\undefined%
\gdef\SetFigFont#1#2#3#4#5{%
  \reset@font\fontsize{#1}{#2pt}%
  \fontfamily{#3}\fontseries{#4}\fontshape{#5}%
  \selectfont}%
\fi\endgroup%
\begin{picture}(8550,8164)(2326,-8519)
\put(2326,-2311){\makebox(0,0)[lb]{\smash{\SetFigFont{5}{6.0}{\rmdefault}{\mddefault}{\updefault}{\color[rgb]{0,0,0}$P_{i+1}^++Q_{j+1}^0$}%
}}}
\put(5326,-511){\makebox(0,0)[lb]{\smash{\SetFigFont{5}{6.0}{\rmdefault}{\mddefault}{\updefault}{\color[rgb]{0,0,0}$P_{i+1}^++Q_{j+2}^0$}%
}}}
\put(10801,-661){\makebox(0,0)[lb]{\smash{\SetFigFont{5}{6.0}{\rmdefault}{\mddefault}{\updefault}{\color[rgb]{0,0,0}$P_{i+2}^++Q_{j+2}^0$}%
}}}
\put(9226,-2461){\makebox(0,0)[lb]{\smash{\SetFigFont{5}{6.0}{\rmdefault}{\mddefault}{\updefault}{\color[rgb]{0,0,0}$P_{i+2}^++Q_{j+1}^0$}%
}}}
\put(4351,-8386){\makebox(0,0)[lb]{\smash{\SetFigFont{5}{6.0}{\rmdefault}{\mddefault}{\updefault}{\color[rgb]{0,0,0}$P_{i+1}^0+Q_{j+1}^-$}%
}}}
\put(10876,-6586){\makebox(0,0)[lb]{\smash{\SetFigFont{5}{6.0}{\rmdefault}{\mddefault}{\updefault}{\color[rgb]{0,0,0}$P_{i+2}^0+Q_{j+2}^-$}%
}}}
\put(9151,-8461){\makebox(0,0)[lb]{\smash{\SetFigFont{5}{6.0}{\rmdefault}{\mddefault}{\updefault}{\color[rgb]{0,0,0}$P_{i+2}^0+Q_{j+1}^-$}%
}}}
\put(7276,-6361){\makebox(0,0)[lb]{\smash{\SetFigFont{5}{6.0}{\rmdefault}{\mddefault}{\updefault}{\color[rgb]{0,0,0}$P_{i+1}^0+Q_{j+2}^-$}%
}}}
\end{picture}

%% file: sigmaijdisc.pstex_t
\begin{picture}(0,0)%
\includegraphics{sigmaijdisc.pstex}%
\end{picture}%
\setlength{\unitlength}{1579sp}%
\begingroup\makeatletter\ifx\SetFigFont\undefined%
\gdef\SetFigFont#1#2#3#4#5{%
  \reset@font\fontsize{#1}{#2pt}%
  \fontfamily{#3}\fontseries{#4}\fontshape{#5}%
  \selectfont}%
\fi\endgroup%
\begin{picture}(8550,8164)(2326,-8519)
\put(2326,-2311){\makebox(0,0)[lb]{\smash{\SetFigFont{5}{6.0}{\rmdefault}{\mddefault}{\updefault}{\color[rgb]{0,0,0}$P_{i+1}^++Q_{j+1}^0$}%
}}}
\put(5326,-511){\makebox(0,0)[lb]{\smash{\SetFigFont{5}{6.0}{\rmdefault}{\mddefault}{\updefault}{\color[rgb]{0,0,0}$P_{i+1}^++Q_{j+2}^0$}%
}}}
\put(10801,-661){\makebox(0,0)[lb]{\smash{\SetFigFont{5}{6.0}{\rmdefault}{\mddefault}{\updefault}{\color[rgb]{0,0,0}$P_{i+2}^++Q_{j+2}^0$}%
}}}
\put(9226,-2461){\makebox(0,0)[lb]{\smash{\SetFigFont{5}{6.0}{\rmdefault}{\mddefault}{\updefault}{\color[rgb]{0,0,0}$P_{i+2}^++Q_{j+1}^0$}%
}}}
\put(4351,-8386){\makebox(0,0)[lb]{\smash{\SetFigFont{5}{6.0}{\rmdefault}{\mddefault}{\updefault}{\color[rgb]{0,0,0}$P_{i+1}^0+Q_{j+1}^-$}%
}}}
\put(10876,-6586){\makebox(0,0)[lb]{\smash{\SetFigFont{5}{6.0}{\rmdefault}{\mddefault}{\updefault}{\color[rgb]{0,0,0}$P_{i+2}^0+Q_{j+2}^-$}%
}}}
\put(9151,-8461){\makebox(0,0)[lb]{\smash{\SetFigFont{5}{6.0}{\rmdefault}{\mddefault}{\updefault}{\color[rgb]{0,0,0}$P_{i+2}^0+Q_{j+1}^-$}%
}}}
\put(7276,-6361){\makebox(0,0)[lb]{\smash{\SetFigFont{5}{6.0}{\rmdefault}{\mddefault}{\updefault}{\color[rgb]{0,0,0}$P_{i+1}^0+Q_{j+2}^-$}%
}}}
\end{picture}